\documentstyle[10pt,a4,amscd,amssymb,verbatim]{amsart} 
\newtheorem{Thm}{Theorem}

\newtheorem{lemma}{Lemma}

\newcommand{\Par}{\operatorname{Par}}
\newcommand{\LGZ}{\operatorname{LGZ}}
\newcommand{\coeff}{\operatorname{coeff}}

\newcommand{\Hom}{\operatorname{Hom}}
\newcommand{\End}{\operatorname{End}}

\newcommand{\modfg}{\operatorname{-modfg}}
\newcommand{\modf}{\operatorname{-modfr}}

\newcommand{\Span}{\operatorname{span}}

\newcommand{\Ind}{\operatorname{Ind}}

\newcommand{\rad}{\operatorname{rad}}
\newcommand{\GZ}{\operatorname{GZ}_n}
\newcommand{\GZQ}{\operatorname{GZ}_{\Q}}

\newcommand{\rank}{\operatorname{rank}}
\newcommand{\Std}{\operatorname{Std}}
\newcommand{\T}{\operatorname{Tor}}
\newcommand{\Fr}{\operatorname{Fr}}
\newcommand{\RStd}{\operatorname{RStd}}
\newcommand{\Class}{\mathfrak C}

\newcommand{\Z}{{\mathbb Z}}
\newcommand{\Q}{{\mathbb Q}}
\newcommand{\A}{ R S_n}
\newcommand{\Ad}{ R S^{\,\ast}_n}
\newcommand{\AQ}{\Q S_n}
\newcommand{\AF}{\F S_n}
\newcommand{\F}{{\mathbb F}_p}

\newcommand{\Aup}{ {R} S_n^{\, \rhd \lambda}}
\newcommand{\Aupe}{ {R} S_n^{\, \unrhd \lambda}}
\newcommand{\AQup}{ {\mathbb Q} S_n^{ \,  \rhd \lambda}}
\newcommand{\AQupe}{ {\mathbb Q} S_n^{ \,  \unrhd \lambda}}
\newcommand{\AFup}{ {\mathbb F}_p S_n^{ \, \rhd \lambda}}
\newcommand{\AFupe}{ {\mathbb F}_p S_n^{ \, \unrhd \lambda}}

\edef\savecatcodeat{\the\catcode`@}
\catcode`\@=11

\def\tb@ifSpecChars#1#2{#1}
\def\tb@ifNoSpecChars#1#2{#2}

\def\tableau{%
  \bgroup
  \@ifstar{\let\Tif\tb@ifNoSpecChars\tb@tableauB}
          {\let\Tif\tb@ifSpecChars\tb@tableauB}}

\def\tb@tableauB{
  \@ifnextchar[{\tb@tableauC}{\tb@tableauC[]}}

\def\tb@tableauC[#1]{\hbox\bgroup%
    \let\\=\cr
    \def\bl{\global\let\tbcellF\tb@cellNF}%
    \def\tf{\global\let\tbcellF\tb@cellH}
    \dimen2=\ht\strutbox \advance\dimen2 by\dp\strutbox%
    \ifx\baselinestretch\undefined\relax%
    \else%
       \dimen0=100sp \dimen0=\baselinestretch\dimen0%
       \dimen2=100\dimen2 \divide\dimen2 by\dimen0%
    \fi%
    \let\tpos\tb@vcenter
    \tb@initYoung
    \tb@options#1\eoo
    \let\arrow\tb@arrow%
    \dimen0=\Tscale\dimen2%
    \dimen1=\dimen0 \advance\dimen1 by \tb@fframe%
    \lineskip=0pt\baselineskip=0pt
    \def\tb@nothing{}%
    \def\endcellno{$\rss\egroup\bss\egroup}
    \def\endcell{\endcellno\kern-\dimen0}
    \def\begincell{\vbox to\dimen0\bgroup\vss\hbox to\dimen0\bgroup\hss$}%
    \let\overlay\tb@overlay%
    \let\fl\tb@fl%
    \let\lss\hss\let\rss\hss\let\tss\vss\let\bss\vss
    \def\mkcell##1{
        \let\tbcellF\tb@cellD
        \def\tb@cellarg{##1}
        \ifx\tb@cellarg\tb@nothing\let\tb@cellarg\tb@cellE\fi%
        \begincell\tb@cellarg\endcellno
        \tbcellF}
    \let\savecellF\tbcellF
     \Tif{\catcode`,=4\catcode`|=\active}{}\tb@tableauD}%

\let\tb@savetableauD\tableauD
{
    \catcode`|=\active \catcode`*=\active \catcode`~=\active%
    \catcode`@=\active
\gdef\tableauD#1{%
  \Tif{
    \mathcode`|="8000 \mathcode`*="8000%
    \mathcode`~="8000 \mathcode`@="8000%
    \def@{\bullet}%
    \let|\cr
    \let*\tf
    \let~\sk
  }{}%
  \tpos{\tabskip=0pt\halign{&\mkcell{##}\cr#1\crcr}}%
  \global\let\tbcellF\savecellF
  \egroup
  \egroup}
}
\let\tb@tableauD\tableauD
\let\tableauD\tb@savetableauD
\let\tb@savetableauD\undefined

\def\tb@options#1{\ifx#1\eoo\relax\else\tb@option#1\expandafter\tb@options\fi}

\def\tb@option#1{%
  \if#1t\let\tpos\tb@vtop\fi
  \if#1c\let\tpos\tb@vcenter\fi
  \if#1b\let\tpos\vbox\fi
  \if#1F\tb@initFerrers\fi
  \if#1Y\tb@initYoung\fi
  \if#1s\tb@initSmall\fi
  \if#1m\tb@initMedium\fi
  \if#1l\tb@initLarge\fi
  \if#1p\tb@initPartition\fi
  \if#1a\tb@initArrow\fi
}

\def\tb@vcenter#1{\ifmmode\vcenter{#1}\else$\vcenter{#1}$\fi}

\def\tb@vtop#1{\hbox{\raise\ht\strutbox\hbox{\lower\dimen0\vtop{#1}}}}

\def\tb@initPartition{\def\Tscale{.3}}
\def\tb@initSmall{\def\Tscale{1}}
\def\tb@initMedium{\def\Tscale{2}}
\def\tb@initLarge{\def\Tscale{3}}

\def\tb@initArrow{\dimen2=1.25em}

\def\tb@initYoung{%
  \def\tb@cellE{}
  \let\tb@cellD\tb@cellN
  \def\sk{\global\let\tbcellF\tb@cellNF}}
\def\tb@initFerrers{%
  \def\tb@cellE{\bullet}
  \let\tb@cellD\tb@cellNF
  \def\sk{\bullet}}

\tb@initMedium

\def\tb@sframe#1{%
  \vbox to0pt{
    \vss
    \hbox to0pt{%
      \hss
      \vbox to\dimen1{
        \hrule depth #1 height0pt
        \vss
        \hbox to\dimen1{
          \vrule width #1 height\dimen1
          \hss
          \vrule width #1
          }%
        \vss
        \hrule height #1 depth 0in
        }%
      \kern-\tb@hframe
      }%
    \kern-\tb@hframe}}

\def\tb@hframe{.2pt}\def\tb@fframe{.4pt}\def\tb@bframe{1.2pt}
\def\tb@cellH{\tb@sframe{\tb@bframe}}       
\def\tb@cellNF{}                            
\def\tb@cellN{\tb@sframe{\tb@fframe}}       
\let\tbcellF\tb@cellN                       

\def\tb@rpad{1pt}
\def\tb@lpad{1pt}
\def\tb@tpad{1.8pt}
\def\tb@bpad{1.8pt}

\def\tb@overlay{\endcell\@ifnextchar[{\tb@overlaya}{\begincell}}
\def\tb@overlaya[#1]{\vbox to\dimen0\bgroup%
  \tb@overlayoptions#1\eoo%
  \tss\hbox to\dimen0\bgroup\lss}
\def\tb@overlayoptions#1{\ifx#1\eoo\relax\else\tb@overlayoption#1\expandafter\tb@overlayoptions\fi}

\def\tb@overlayoption#1{
  \if#1t\def\tss{\vskip\tb@tpad}\let\bss\vss\fi
  \if#1c\let\tss\vss\let\bss\vss\fi
  \if#1b\def\bss{\vskip\tb@bpad}\let\tss\vss\fi
  \if#1l\def\lss{\hskip\tb@lpad}\let\rss\hss\fi
  \if#1m\let\lss\hss\let\rss\hss\fi
  \if#1r\def\rss{\hskip\tb@rpad}\let\lss\hss\fi
}

\def\tb@fl{\endcell\begincell\vrule depth 0pt width \dimen0 height \dimen0 \endcell\begincell}

\@ifundefined{diagram}{}{
\def\tb@arrowpad{.5}

\newoptcommand{\tb@arrow}{\@ne}[2]{%
  \endcell
   \begingroup%
   \let\dg@getnodesize\tb@getnodesize
   \dg@USERSIZE=#1\relax%
   \ifnum\dg@USERSIZE<\@ne \dg@USERSIZE=\@ne \fi%
   \dg@parse{#2}%
   \dg@label{\tb@draw{#1}{#2}}}

\def\tb@getnodesize#1#2#3#4#5{\dimen3=\tb@arrowpad\dimen2 #4=\dimen3 #5=\dimen3\relax}
\def\tb@getnodesize#1#2#3#4#5{\ifnum#2=0\ifnum#3=0\tb@getnodesizetail{#4}{#5}\else\tb@getnodesizehead{#4}{#5}\fi\else\tb@getnodesizehead{#4}{#5}\fi}
\def\tb@getnodesizetail#1#2{\dimen3=.5\dimen2 #1=\dimen3 #2=\dimen3}
\def\tb@getnodesizehead#1#2{\dimen3=.5\dimen2 #1=\dimen3 #2=\dimen3}

\def\tb@draw#1#2#3#4{%
        \dg@X=0\dg@Y=0\dg@XGRID=1\dg@YGRID=1\unitlength=.001\dimen0%
        \dg@LBLOFF=\dgLABELOFFSET \divide\dg@LBLOFF\unitlength%
        \dg@drawcalc
        \begincell
        \let\lams@arrow\tb@lams@arrow
        \begin{picture}(0,0)\begingroup\dg@draw{#1}{#2}{#3}{#4}\end{picture}%
        \endcell
        \endgroup
        \begincell}
}

\def\tb@lams@arrow#1#2{%
 \lams@firstx\z@\lams@firsty\z@
 \lams@lastx#1\relax\lams@lasty#2\relax
 \lams@center\z@
 \N@false\E@false\H@false\V@false
 \ifdim\lams@lastx>\z@\E@true\fi
 \ifdim\lams@lastx=\z@\V@true\fi
 \ifdim\lams@lasty>\z@\N@true\fi
 \ifdim\lams@lasty=\z@\H@true\fi
 \NESW@false
 \ifN@\ifE@\NESW@true\fi\else\ifE@\else\NESW@true\fi\fi
 \ifH@\else\ifV@\else
  \lams@slope
  \ifnum\lams@tani>\lams@tanii
   \lams@ht\ten@\p@\lams@wd\ten@\p@
   \multiply\lams@wd\lams@tanii\divide\lams@wd\lams@tani
  \else
   \lams@wd\ten@\p@\lams@ht\ten@\p@
   \divide\lams@ht\lams@tanii\multiply\lams@ht\lams@tani
  \fi
 \fi\fi
 \ifH@  \lams@harrow
 \else\ifV@ \lams@varrow
 \else \lams@darrow
 \fi\fi
}

\catcode`\@=\savecatcodeat
\let\savecatcodeat\undefined

\begin{document}    
\large

\title{Young's seminormal form and simple 
modules for $S_n$  in characteristic $p$.} 
\author{Steen Ryom-Hansen}
\address{Instituto de Matem\'atica y F\'isica, Universidad de Talca \\
Chile\\ steen@@inst-mat.utalca.cl }
\footnote{Supported in part by Programa Reticulados y Simetr\'ia, by FONDECYT grant 1090701 and by 
the MathAmSud project OPECSHA 01-math-10.} 
\maketitle

\begin{abstract}
We realize the integral Specht modules for the symmetric group $ S_n $ 
as induced modules from the subalgebra of the group algebra 
generated by the Jucys-Murphy elements.  
We deduce from this that the simple modules for $ \F S_n $ 
are generated by reductions modulo $ p $ of the corresponding Jucys-Murphy idempotents. 
\end{abstract}

\section{Introduction.} 
This article is a continuation of the investigation 
pursued in [RH1-2] that seeks to demonstrate the importance of Young's seminormal basis
for the modular, that is characteristic $ p $, representation theory of the symmetric group
$ S_n$.
A main obstacle is here 
that Young's seminormal basis is defined over the field $ \Q $ 	
and indeed 
there seems to be a general consensus that 
this obstacle makes 
Young's seminormal basis a characteristic zero phenomenon, essentially. Still
we believe that 
Young's seminormal basis is a fundamental object 
for the modular representation theory as well, and 
we think that the results of our works provide strong evidence in favor of this claim.

\medskip
Let $ \Par_n $ be 
the set of partitions of $n$ and let $ S(\lambda ) $ be 
the integral Specht module for $ S_n $ associated with $ \lambda \in \Par_n $.
Then, as has been known for a long time, 
the set of $ S_{\Q} (\lambda) := S(\lambda ) \otimes_{\Z} \Q $ classifies the 
irreducible $ \Q S_n $-modules when $ \lambda \in \Par_n $, whereas 
the reduced Specht modules $ S(\lambda ) \otimes_{\Z} \F $ are
reducible in general. In fact the irreducible 
modules $ D(\lambda) $ 
for $ \F S_n $ are classified by the set of $ p$-regular partitions $ \Par_n^{reg} $ and
are obtained as 
$ D(\lambda ) = S(\lambda)/ \rad (\cdot, \cdot) $ where $ (\cdot, \cdot)$ is a certain symmetric
bilinear and $ S_n $-invariant form on $S(\lambda) $. 

\medskip
The decomposition numbers $ [ S(\lambda): D(\mu)] $ for $ \F S_n  $ have been the topic
of much research activity in recent years, but still remain unknown in general and even 
the dimensions of $ D(\mu) $ are not known in general. But using 
the theory of Young's seminormal form we obtain in this work, as our 
main Theorem {\ref{simples}}, a 
construction of $ D(\mu) $  
that may be a good starting point for obtaining 
combinatorial expressions for $ \dim D(\mu) $.

\medskip
The 
basic principles behind this construction are 
parallels of standard methods in the modular representation theory of 
algebraic groups. Indeed, let $ S(\lambda)^{ \circledast} $ denote the contragredient dual 
of $ S(\lambda) $. Then $ (\cdot, \cdot)$ corresponds to a homomorphism
$ c_{\lambda}: S(\lambda) \rightarrow  S(\lambda)^{ \circledast} $.
Moreover, for $ \lambda \in \Par_n^{reg}$  we have that 
$ D(\lambda) = im \, \overline{ \,c_{\lambda} } $ where $ \overline{\,c_{\lambda} }$ is the 
reduced homomorphism modulo $ p$. Passing to the representation theory 
of an algebraic group $ G $ over an algebraically closed field of characteristic $p$,  
the Weyl module $ \Delta(\lambda) $, the dual Weyl module $ \nabla(\lambda) $ and 
the simple module $ L(\lambda) $ 
correspond to $ S(\lambda) $, $ S(\lambda)^{ \circledast} $ and 
$ D(\lambda) $ and the bilinear form $ (\cdot, \cdot)$ on $ S(\lambda) $ corresponds 
to a form on $ \Delta(\lambda) $ that we denote the same way. It induces a $ G $-linear
homomorphism $ c_{\lambda}: \Delta(\lambda) \rightarrow \nabla(\lambda) $ 
and the simple module satisfies 
$ L(\lambda) =  im  \, c_{\lambda} $. 
But in the $ G $-module setting, 
$ \nabla(\lambda) $ can also be constructed 
as the module of global sections of a line bundle on the associated flag manifold, 
and using this, one obtains a new 
construction of $ c_{\lambda}$ without using the bilinear form.
The properties of this  
new construction of $ c_{\lambda}$ then provide 
a useful method for obtaining information on $ L(\lambda) $, see
eg. [A, Jan].

\medskip
Returning to the symmetric group, we then look for a
different construction of $ c_{\lambda} $.
For this we prove in our Theorem \ref{promised} that $ S(\lambda) $ 
is induced from a certain subalgebra, denoted $ \GZ$, of the group algebra, 
corresponding 
to the fact that 
$ \nabla(\lambda) $ is induced from a Borel subgroup of $ G$.
Given this, our 
new construction of $ c_{\lambda}$ is obtained from 
a Frobenius reciprocity argument.

\medskip
At the basis of our work are 	
the famous Jucys-Murphy elements $ L_k, \,k = 1,2, \ldots, n $ 
that 
were introduced independently by Jucys and Murphy in [Ju1-3] and [Mu81]. 
They give rise to 
idempotents $ E_t $ of $ \Q S_n $, the Jucys-Murphy idempotents, indexed by  
$ \lambda $-tableaux $ t$,
that are closely related to Young's seminormal basis of the Specht module 
$ S_{\Q}(\lambda) $.
Moreover 
they commute with each other and therefore generate    
a commutative subalgebra of the group algebra. This is the 
algebra $ \GZ$ that was mentioned above, the Gelfand-Zetlin algebra.
In the case of the ground field $ \Q $
it was considered by Okounkov and 
Vershik in [OV] as a kind of  
Cartan subalgebra of a semisimple Lie algebra, but for us it is
important to work with an integral version of $ \GZ $, where the analysis 
of [OV] fails.

\medskip
We have now formulated the main ingredients of our result. 
The surprisingly simple final result is that
$ D(\lambda) $ is generated by $ a_{\lambda}E_{\lambda} $ 	
where $ E_{\lambda} = E_{t^{\lambda}}$ and
$ a_{\lambda} $ is the least common multiple of the denominators of 	
$ E_{\lambda}$. It should be noted that, even though 
it appears to be a very natural idea to 
investigate the $ \Z S_n $ or $ \F S_n $-submodule
of $ S(\lambda) $ generated by $ a_{\lambda}E_{\lambda} $, 
the only reference in the literature along these lines is [RH2], as far as 
we know.

\medskip	
In an important recent paper [BK], J. Brundan and A. Kleshchev showed that 
$ \F S_n$ is a $ \Z$-graded algebra 
in a nontrivial way 
by establishing an isomorphism between 
$ \F S_n$
and the cyclotomic KLR-algebra, i.e. cyclotomic Khovanov-Lauda-Rouquier 
algebra, of type $ A$.
Their results work in greater
generality than $ \F S_n$ but we shall only consider this case.
J. Hu and A. Mathas refined in [HuMa] this graded structure on $ \F S_n$ to  
a graded cellular algebra structure by constructing an explicit graded cellular basis. 
A second goal of our paper is to 
show that key features of their constructions can be 
carried out entirely within the theory of Young's seminormal form, as developed by 
Murphy.
We hope that this approach to their results, together with
our main Theorem {\ref{simples}}, may provide a combinatorial expression for 
$ \dim D(\lambda) $.

\medskip
The generators of the cyclotomic KLR-algebra are 
$$ \{ \, e( {\bf i}) \,| \,{\bf i} \in (\F)^n \} \cup \{ y_1, \ldots, y_{n} \}
\cup \{ \psi_1, \ldots , \psi_{n-1} \}	$$
and [BK] prove their Theorem by constructing elements in $ \F S_n $, 
denoted the same way, that 
verify the cyclotomic KLR-algebra relations.
The $ y_i $ are essentially Jucys-Murphy operators and 
$  e( {\bf i}) $ are certain idempotents, not necessarily nonzero. In fact they can be identified with the
idempotents 
constructed in [Mu83] 
by summing Jucys-Murphy idempotents $ E_t $
over tableaux classes. 
The elements $ \psi_i $ are the most difficult to handle 
and [BF] take as starting point for this 
certain explicitly given intertwining elements $ \phi_i $. 
These intertwiners, 
together with the $ e( {\bf i}) $ and $ y_i $, 
already satisfy relations that are close to the 
cyclotomic KLR-relations but still need to be adjusted to get the complete match.

\medskip
We here give a natural construction of the intertwining elements $ \phi_i $ within 
Murphy's theory for the seminormal basis. 
Indeed, we see them as natural analogues of certain 
elements $ \Psi_i  $ of the Hecke algebra that 
appear in [Mu92], 
although only in the semisimple case.
We show that Murphy's ideas, in a suitable sense, 
can be carried out over $ \F $ as well.
From this we obtain a cellular basis for $ \F S_n$ 
using a 
modification of the construct ions done in 
[HuMa].	 

\medskip
Let us sketch the layout of the paper. In section 2
we fix the basic notation of the paper. It is mostly standard, except possibly for 
the notion of tableau class which was introduced in [Mu83]. We also review
the construction from [Mu83] of the tableau class idempotents.
Section 3 contains the construction of the intertwiners $ \Psi_{L,i} $.
This requires 
a control of the denominators of the Jucys-Murphy idempotents that are involved in 
the tableau class idempotents. In section 4 we construct the cellular basis. 
In section 5 we first introduce the Gelfand-Zetlin algebra $ \GZ $ and then set up the induction 
functor. We then prove that the Specht module is induced up from a ``rank one'' module 
of $ \GZ $. An important ingredient for this is a uniqueness statement, due to James [J], 
of the integral Specht module. 
Finally in section 6 we deduce our main results.

\medskip
Note that the notation used throughout the paper may 
vary slightly from the one used in the introduction.

\medskip
It is a pleasure to thank H. H. Andersen, J. Brundan, P. Desrosier, S. Griffeth, L. Lapointe, 
A. Mathas, O. Mathieu, D. Plaza and 
W. Soergel, among others, for useful discussions. Finally, it is a special pleasure to thank the referee for 
his/her useful suggestions and for pointing out several inaccurracies in a first version of the paper.


\section{Basic notation and idempotents in positive characteristic.}
We are concerned with the representation theory of the symmetric group 
$ S_n $ in positive characteristic. 
Let us first set up the basic notation. 
Let $ p$ be a prime. We use the ground rings $ R $, $ \Q $ and $  \F  $, 
where $ R $ is the localization of $ \Z $ at the prime ideal $ (p) $ and where $  \F  $ is
the finite field of $ p $ elements. 
Then $ R $ is a local ring with maximal ideal $ p R $ and
$ R/pR = { \F} $.
Let $ n $ be a positive integer and let $ S_n $ be the symmetric 
group on $n$ letters. 
An $n$-composition
is a sequence $ \lambda = (\lambda_1 , \lambda_2, \ldots , \lambda_k) $ 
of positive integers with sum $ n$. 
An $ n $-partition is an $ n$-composition 
$ \lambda = (\lambda_1 , \lambda_2, \ldots , \lambda_k) $ 
such that $ \lambda_i \geq  \lambda_{i+1} $ for all $ i $.
The set of $ n $-partitions is denoted $ \Par_n$. For $ \lambda \in \Par_n$,  
the associated  
Young diagram, also denoted $ \lambda$, is the graphical representation 
of $ \lambda $ through $ n $ empty boxes in 
the plane. The first $ \lambda_1 $ boxes are placed in the 
first row, the next $ \lambda_2 $ boxes are placed in the
second row, left aligned with respect to the first row, etc. 
This is the 
English notation for Young diagrams. 
The boxes are denoted the {\it nodes} of $ \lambda $ and 
are indexed using matrix convention. Thus the node of $ \lambda $ 
indexed by $ [2,3] $ is the one 
situated in the second row and the third column of $ \lambda$. 
The $ p $-residue diagram of a partition $ \lambda $ is obtained 
by writing $ j-i  \mbox{ mod } p $ in the $ [i, j] $'th node of $ \lambda $.
The $ [i, j] $'th node is called a $ k $-node of $ \lambda $ if $ k = j -i \mbox{ mod } p $.



\medskip
Let $ t $ be a $ \lambda$-tableau, i.e. a filling 
of the nodes of $ \lambda $ using the numbers of $ \{1, 2, \ldots, n \} $, each once.
We write 
$ t[i,j] = k $ if the $ [i,j] $'th node of $ t $ is filled in with $ k$
and $ r_{ t}( k ) = j-i $ if $ t[i, j ] = k$. Then $ r_{ t}( k )$ is also referred to as the content of $ t $ 
at the node containing $k$. For $ k \in \{1, 2, \ldots, n \} $ 
we define $ t(k) := [i,j] $ where $ t[i,j] = k $. 
A tableau $ t $ is said to be row standard if $ t[i, j ] \leq t[i, j+1 ] $ 
for all $ i, j $ such that the terms are defined.
The set of row standard tableaux of all $ n $-partitions is denoted $ \RStd(n) $.
A tableau $ t $ is said to be standard if $ t[i, j ] \leq t[i, j+1 ] $ and $ t[i, j ] \leq t[i+1, j ] $
for all $ i, j $ such that the terms are defined.
For $ \lambda $ an $n$-partition, we 
define $ Shape(t) := \lambda $ if $ t $ is a $ \lambda$-tableau.
The set of standard tableaux of $ n $-partitions is denoted $ \Std(n) $ and the set of 
standard tableaux of shape $ \lambda $ is denoted $ \Std(\lambda)$.

\medskip
If $ t $ and $ s $ are tableaux we write  
$ t \sim_p s $ if $ r_t(k) = r_s(k) \mbox{ mod} \, p $ whenever $ t[i,j ] = s[i_1, j_1] = k $. 
This defines an equivalence relation on the set of all tableaux which can be restricted to an 
equivalence relation on the set of standard tableaux. 
When we refer to a tableau class we always mean a class with 
respect to the last relation, consisting of standard tableaux. 
If $ t $ is a standard tableau we denote its class by $ [t] $.
In general we refer to tableau classes using capital letters, like $ T $ or $S $. 
We denote by $ {\Class}_n $ the set of tableau classes of all 
$ n$-partitions. 

\medskip
For $ t $ a tableau, 
the multiset of residues $ \{  r_t(k) \,  |\, k=1, \ldots, n \} $ depends only on $ Shape(t) $ and this
induces an equivalence relation on $ \Par_n $ that we also denote $ \sim_p$. 
The blocks of $ {\mathbb F} S_n $ are given by this equivalence relation according to the Nakayama conjecture, 
see for example Wildon's notes [W]. 
Clearly, for two tableaux $ s, t$ we have that $ t \sim_p s $ implies $ Shape(t) \sim_p Shape(s)$.

\medskip
We use the convention that $ S_n $ acts on the right on $ \{1, \ldots, n \} $ and
hence on tableaux. 
In other words, we multiply 
cycles in $ S_n $ from the left to the right.

\medskip
For $ t $ a $ \lambda$-tableau, we define the associated element $ d(t) \in S_n $ by
$$  t^{\lambda} d(t)= t$$
where 
$ t^{ \lambda} $ denotes the highest $ \lambda $-tableau, having the numbers
$ \{ 1,2, \ldots, n \} $ filled in along rows. Highest refers to the dominance order 
$ \unlhd $ 
on tableaux. It is derived from the dominance order $ \unlhd $ on compositions
given by
$$  \lambda   \unlhd \mu   \mbox{ if }
\sum_{i=1}^{m} \lambda_i \leq \sum_{i=1}^{m} \mu_i \, \, \mbox{   for } m = 1, 2, \ldots ,
\mbox{min}(k, l) $$
for $ \lambda = (\lambda_1, \ldots, \lambda_k) $ and 
$ \mu = (\mu_1, \ldots, \mu_l) $ 
by 
viewing tableaux as series of compositions.
Similarly the dominance order can be extended to pairs of tableaux 
in the following way
$$ (s, t) \unlhd   (s_1, t_1 ) \mbox{ if } s \unlhd  s_1 \mbox{ and } t \unlhd t_1. $$
In [HuMa] this order on pairs of tableaux is called the strong dominance order and is written
$ \blacktriangleleft $. 
The dominance orders are all partial. 

\medskip
Let $ t $ be a $ \lambda $-tableau with node $ (i,j) $. 
The $ (i, j)$-hook consists 
of the nodes to the right and below the $ (i, j ) $ node, its 
cardinality is the hook-length $ h_{i,j} $. The product of all hook-lengths 
only depends on $ \lambda $ and is denoted $ h_{\lambda} $.
The hook-quotient is 
$ \gamma_{t,n} = \prod  \frac{h_{i,j}}{h_{i,j}-1} $ 
with the product taken over 
all nodes in the row of $ \lambda $ that contains $ n $, omitting hooks of length one. 
For general $ i $, we define $ \gamma_{t,i} $ similarly, by first deleting from $ t $ the nodes containing 
$ i+1, i+2, \ldots , n $. 
We set $ \gamma_t = \prod_{i=2}^n \gamma_{t,i} $. 

\medskip 
In general, when we use $ \lambda $ as a subscript it refers to the tableau $ t^{\lambda} $.
In this situation we have
$$ \gamma_{\lambda } = \gamma_{t^{\lambda} } = \prod_i \lambda_i !. $$

For $ k = 1,2, \ldots, n $ the Jucys-Murphy elements $ L_k \in {\mathbb	 Z} S_n $ are
defined by
$$ L_k := (1, k ) + (2,k ) + \ldots + (k-1, k) $$
with the convention that $ L_1  := 0 $. 
They commute with each other and satisfy the following commutation relations with the 
simple transpositions
\begin{equation}{\label{Jucys-Murphy}}
\begin{array}{lr}
(k-1,k) L_k  = L_{k-1} (k-1,k) +1 \\
(k-1,k) L_{k-1}  = L_{k} (k-1,k) -1 \\
(k-1,k) L_l  = L_l (k-1,k) & \mbox{     if } l \not= k-1, l \not= k. \\
\end{array}
\end{equation}
These elements are a key ingredient for understanding the 
representation theory of $ S_n$. Their generalizations appear in many contexts of representation theory, for 
example as 
the degenerate affine Hecke algebra, where the $ L_k $ 
are commuting generators that satisfy the above relations with the simple transpositions.
In the original works of Jucys and 
Murphy, [Ju1], [Ju2], [Ju3] and [Mu81], the $ L_k $'s were used to construct orthogonal 
idempotents 
$ E_t \in { \mathbb Q} S_n $, indexed by tableaux $ t $, and to derive Young's seminormal form 
from them.
We denote these idempotents the Jucys-Murphy idempotents.
Their construction is as follows 
$$ E_t := \prod_{  \{c \, |    -n \, < \,  c \,  < n \} }  \prod_{  \{ \, i \, |  r_{ t}(i) \not= c  \} \,} 
\frac{ L_i -c }{ r_{t}(i) -c}.$$
For $ t $ standard we have $ E_t  \not= 0 $, whereas 
for $ t $ nonstandard either $ E_t = 0 $, or 
$ E_t = E_s $ for some standard tableau $ s $ related to $ t$, see [Mu83] page 260. 
For example, if $ t $ is obtained from the standard tableau $ s $ by interchanging $ k-1 $ and $ k $
that occur in the same row or column of $ s $ then $ E_t =0 $.
Running over all standard tableaux, the $ E_t $ form a set of primitive and 
complete idempotents, that is their sum is $1$.
Moreover, they are eigenvectors for the action of the Jucys-Murphy operators in $ {\mathbb Q} S_n $,
since 
\begin{equation}{\label{Mu_page_506}}
(L_k-r_t(k)) E_t =0 \mbox{ or equivalently }
L_k = \sum_{ t \in \Std(n)} r_t(k) E_t \end{equation}
which is the key formula for deriving Young's seminormal basis  
from them. In this situation (\ref{Jucys-Murphy}) gives Young's seminormal form 
for the action of $ \sigma_i $ on the seminormal basis. 

\medskip
Unfortunately, the $ E_t $ contain many 
denominators and hence it is not possible to reduce them modulo $ p $. 
In order to overcome this obstacle, Murphy introduced in [Mu83] certain 
elements $ E_T $ for each tableau class $ T $.
They are defined as follows
\begin{equation}{\label{ET}} E_T := \sum_{ t \in T } E_t. \end{equation}
He showed that the $ E_T $'s, with $ T $ varying over all classes, give a set of complete
orthogonal idempotents in $ \A $.
The most difficult part of this is to show that 
$ E_T \in \A $ since they are clearly orthogonal, idempotent and complete (note at this point that in [Mu83] it is 
not stated clearly that 
the sum in (\ref{ET}), going over the elements of the tableau class $T$, should only involve standard tableaux). 
We now present his proof that $ E_T \in \A $, in our notation. Several of 
its ingredients 
will be important for us. See also [MaSo] for a presentation of this and related results from an abstract point of view.

\medskip
A key 
point is to consider $ F_t $ for $ t $ any tableau, given by 
\begin{equation}{\label{F_t}}
F_t := \prod_{ \{  \, c \, |    -n  < c <  n \} } \prod_{  \{ \, i \, |  r_{ t}(i) \not= c 
\! \! \!
\mod p  \} \,} 
\frac{ L_i -c }{ r_{t}(i) -c}.
\end{equation}
It is clear that $ F_t \in \A $ and that all 
$ F_t $'s and $ E_T$'s commute.
The denominator of $ F_t $ depends only on the underlying 
partition $ Shape(t )= \lambda $ of $ t $ and is denoted $ w^{\lambda} $. Although $ w^{\lambda} $ 
is not constant on the classes, we have that 
$ w^{\lambda} = w^{\mu} \mbox{ modulo } p $ if $ \lambda \sim_p \mu $. Especially, 
if $ s \sim_p t $ we get that $ w^{Shape(s)} = w^{Shape(t)}  \mbox{ modulo } p $.
The numerator of $ F_t $ only depends on the class $ [t] $ of $ t $ and so we have 
$$ F_s = F_t \, \, \, \, \mbox{ if } s \sim_p t   \mbox{ and } Shape(s) = Shape(t) $$
Suppose that $ t \in T $. 
Using (\ref{Mu_page_506}) we get that
\begin{equation}{\label{F_t E_s}}
 F_t E_s = \left\{
\begin{array}{ll}
w^{Shape(s)}/w^{Shape(t)} E_s & \mbox{ if } s \sim_p t \\ 0 & \mbox{ otherwise } 
\end{array} \right.
\end{equation}
and so we deduce 
\begin{equation}{\label{and-so-we deduce}}
F_t = \frac{1}{w^{\lambda}} \sum_{ s \in T} \, w^{Shape(s)}  \, E_s.
\end{equation}
Hence $ E_{T} F_t = F_t $ where we set
$ T= [t] $. 
Using this we get for any positive integer $ m $ that 
$$ 
\begin{array}{c}
(E_T -F_t)^m =  \sum_{i=0}^m 
\left( \begin{array}{c} m \\ i \end{array} \right) 
(-1)^{m-i} E_T^i F_t^{m-i} = \\
E_T - 1 +1 + \sum_{i=0}^{m-1} 
\left( \begin{array}{c} m \\ i \end{array} \right) 
(-1)^{m-i} E_T^i F_t^{m-i} = \\
E_T - 1 +\sum_{i=0}^{m} 
\left( \begin{array}{c} m \\ i \end{array} \right) 
(-1)^{m-i} F_t^{m-i} = E_T-1 + (1-F_t)^m.
\end{array}
$$
Combining this with equation ({\ref{and-so-we deduce}}) we arrive at 
the formula 
\begin{equation}{\label{murphy_formula}}
E_T = 1 - (1-F_t)^m +   \sum_{  s \in T}  \left(1- \frac{w^{Shape(s)} }{w^{\lambda}}\right)^m E_s.
\end{equation}
Using it, the proof that  
$ E_T \in \A $ follows by taking $ m $ big enough for  
$$ (1- \frac{w^{Shape(s)} }{w^{\lambda}})^m E_s \in \A $$ to hold for all $ s \in T $. 

\section{Commutation rules.}
Our first aim is to generalize certain results valid for $ E_t $ to $ E_T$. 
We are especially looking for a generalization 
for $ E_T $ of the elements denoted $ \Psi_t $ in [Mu92]. For this we need 
to work out the commutation relations between $ E_T $ and the simple transpositions 
$ \sigma_k =(k-1,k) $. 

\medskip
Assume that $ t $ and $ t \sigma_k $ are standard tableaux and write 
$ T:= [t] $.
We first consider the case where $ [ t \sigma_k] =  [t] $, that is 
$ r_{t}( k-1 ) = r_{t}(k ) \mbox{ mod } p $. 
We prove the following Lemma. 
\begin{lemma}{\label{commute}} In the above situation $ [\sigma_k t] =  [t] = T $ we have 
$$ \sigma_k E_T = E_T \sigma_k.  $$
\end{lemma}
\begin{pf*}{Proof}                                  
We consider the commutator $ [ \sigma_k, E_T ] $. By the previous section it belongs to 
$\A$. We show that it actually 
belongs to $ p^N \A $ for any
positive (big) integer $ N $, from which the result follows. Fix therefore such an $ N$. 
We use formula ({\ref{murphy_formula}}) and first consider the individual terms of that sum. 
We choose 
$ m $ big enough for $  (1- \frac{ w^{Shape(t_1)}}{{w^{\mu}}}) ^m E_{t_1} \in p^N \A $ 
to hold for all
$ t_1 \in T $. From this we get that 
$$ 
\left[\sigma_k, \sum_{ t_1 \in T}  \left(1- \frac{w^{Shape(t_1)} }{w^{\mu}} \right)^{ \! \! m} \!  E_{t_1} \right ] \in 
p^N \A $$
and so by ({\ref{murphy_formula}}) it is enough to prove that $ \sigma_k $ commutes with 
$ F_t $.  

Now by the commutation rules ({\ref{Jucys-Murphy}}), we have that $ \sigma_k $ commutes
with all terms of $ F_t $ of the form $ L_i - c $ where $ i \not= k-1, k$. 
The remaining terms may be grouped
together in pairs of the form $$ (L_{k-1} -c ) (L_{k} - c) $$ since by assumption
$ r_{t}(k-1 ) = r_{t}(k ) \mbox{ mod } p $. But these expressions are symmetric in $ L_{k-1} $ and $ L_k $ 
and therefore commute with $ \sigma_k $ by the 
commutation rules ({\ref{Jucys-Murphy}}). The Lemma is proved.
\end{pf*}

We next consider the case where $ s, t  $ are both standard and $ s= t \sigma_k   \notin [t] $, that is
$ r_{t}(k-1 ) \not= r_{t} (k ) \mbox{ mod } p $. We set $ S:= [s ] $ and $ T:=[t] $.
In order to work out the commutation rule between $ \sigma $ and $ E_T $ in this case,
we first consider 
$ E := E_S + E_T $. We need the following auxiliary Lemma. 
\begin{lemma}{\label{future_use}} 
$ E $ belongs to $ \A $ and commutes with $ \sigma_k $.
\end{lemma}
\begin{pf*}{Proof}                                  
Clearly $ E $ belongs to $ \A $ since $ E_S $ and $E_T $ do. 
For each $ s \in S $ we have that either $ s \sigma_k \in T $ or $ s \sigma_k $ is nonstandard, and conversely
for each $ t \in T $ we have that either $ t \sigma_k \in S $ or $ t \sigma_k $ is nonstandard. 
Accordingly, the sum $ E = \sum_{t \in T} E_t + \sum_{s \in S} E_s $ may be split into 
terms $ E_s + E_{ s \sigma_k  } $ with $ s \in S $ and $s \sigma_k \in T$, terms $ E_s $ with $  s \sigma_k  $ 
nonstandard and terms $ E_t $ with $  t \sigma_k $ nonstandard. 

The first kind of terms are symmetric in $ L_{k-1} $ and $ L_k$ and 
therefore commutes with $ \sigma_k $. For the second kind of terms, $ k-1 $ and $ k $ occur in the same 
row or column of $ s $, next to each other. Thus, using Theorem 6.4 of [Mu92] 
together with the formula for $ x_{tt} $ of {\it loc. cit.} appearing eight lines above that Theorem 6.4, we get 
that $ E_s \sigma_k = \pm E_s $, where the sign is positive iff $ k-1 $ and $ k $ are in the same row of $ s$.
Applying the antiautomorphism $ \ast$ of $ \A $ that 
fixes the simple transpositions we get $ \sigma_k  E_s = \pm E_s $, with the same 
sign as before, and so indeed $ E_s $ and $ \sigma_k $ commute.
Finally, the third kind of terms involving $ E_s $ are treated the same way. 
\end{pf*} 

\medskip
For each tableau class $ T$
we choose an arbitrary $ t \in T $ and define 
$$ r_T(i) := r_t(i).  $$                                  
Thus $ r_T(i) \in \mathbb Z $, hence also $ r_T(i) \in R$, 
but it is only well defined modulo $ p $. 
With this notation we can formulate our next Lemma.
\begin{lemma}{\label{due_to_Murphy}} 
Assume that $ S $ and $ T $ are chosen as above. Then 
there is a positive integer $ m_1  $ such that the following formulas hold for $ m \geq m_1 $
$$ E_T = \left( \frac{L_k -r_S(k)}{ r_T(k) - r_S(k) } \right)^{ \! m}  \!  E, \, \, \, \, \, \, \, 
E_S = \left( \frac{L_{k-1} -r_S(k)}{ r_T(k) - r_S(k) } \right)^{ \! m} \! E. $$
\end{lemma}
\begin{pf*}{Proof}                                  
For tableaux $ t $ and $ s$ we define $ E_{t, s} = E_t + E_{ s} $. Then 
obviously $ E_{t, s} $ is idempotent and 
by (\ref{Mu_page_506}), we have that 
$ L_k = \sum_{ u \in \Std(n) } r_u(k) E_u $. Hence for standard tableaux $ s $ and $ t $ such 
that $ s = \sigma_k t $ we get that
\begin{equation}{\label{and_similarly_t}}
E_t = \left( \frac{L_{k} -r_s(k)}{ r_t(k) - r_s(k) } \right) E_{t,s}.
\end{equation}
Similarly we have that 
\begin{equation}{\label{and_similarly_s}}
E_s = \left( \frac{L_{k-1} -r_s(k)}{ r_t(k) - r_s(k) } \right) E_{t,s}. 
\end{equation}
Note that ({\ref{and_similarly_t}}) and ({\ref{and_similarly_s}}) hold even 
if only $ s $ or $ t $ is standard, as long as $ s = \sigma_k t $. 
The formulas are used in Murphy's papers but 
unfortunately they do not generalize to $ E_T $ or $ E_S$ 
since we do {\it not} have 
$ L_k = \sum_{ U \in {\Class}_n }\,   r_U(k) E_U $ 
even though $ r_u(k) = r_{u_1}(k) \mbox{ mod } p$ for 
$ u \sim_p u_1 $. The problem is that the individual $ E_u $ lie in $ {\mathbb Q} S_n $ rather 
than $ R  S_n $.

\medskip
In order to solve this problem we proceed as follows. Consider first an $ a \in p R $. 
From the binomial expansion we get the following
formula in $ R[x] $, valid for any positive integer $ m$
\begin{equation}{\label{binomial}}
(x+a)^{p^{m}} = x^{p^{m}} \, \, \mbox{mod} \, \,  p^{m+1}R[x].
\end{equation}
We deduce from it the formula 
\begin{equation}{\label{binomial_fraction}}
\left(\frac{x+a}{c+d}\right)^{p^{m}} = \left(\frac{x}{c}\right)^{p^{m}} 
\, \, \mbox{mod} \,\,  p^{m+1}R[x] 
\end{equation}
for any $ c \in R^{\times} $ and $ a, d \in p R $.

\medskip
Set $ T^{\prime}:= T \cup S s_k $. Then $ T^{\prime} \setminus T $ consists of nonstandard tableaux and 
so $ E_T = \sum_{t \in T^{\prime} } E_t$. 
Using (\ref{and_similarly_t}), for $ m $ is large enough we get
$$ 
\begin{array}{c}
E_T = \sum_{t \in T^{\prime} } E_t = \sum_{t \in T^{\prime}, s= t \sigma_k  } 
\left( \frac{L_{k} -r_s(k)}{ r_t(k) - r_s(k) } \right) E_{t, s} = \\
\sum_{t \in T^{\prime}, s= t \sigma_k  } 
\left( \frac{L_{k} -r_s(k)}{ r_t(k) - r_s(k) } \right)^{p^m} E_{t,s} = 
\sum_{t \in T^{\prime}, s= t \sigma_k  } 
\left( \frac{L_{k} -r_S(k)}{ r_T(k) - r_S(k) } \right)^{p^m} E_{t,s} 
\end{array}
$$
The last equality follows from (\ref{binomial_fraction}), 
since for any $ N$ we may choose $ m $ big enough to make the difference of the two sides 
belong to $ p^N \A $. From the last expression we then get that $ E_T $ is equal to 
$$ 
\begin{array}{c}
\left( \frac{L_{k} -r_S(k)}{ r_T(k) - r_S(k) } \right)^{p^m} \! (E_S + E_T )
=
\left( \frac{L_{k} -r_S(k)}{ r_T(k) - r_S(k) } \right)^{p^m} \! E 
\end{array}
$$
as claimed. The other equality is proved the same way. 
By choosing $ m $ even bigger we obtain an $m_1 $ that works for both equations. 
\end{pf*} 
At this stage Murphy constructs in [Mu92], using 
the formulas ({\ref{and_similarly_t}) and (\ref{and_similarly_s}),  
elements $ \Psi_t $ and $ \Phi_t $ of $ {\mathbb Q} S_n $ 
satisfying 
\begin{equation}{\label{key-property}} E_{\lambda} \Phi_t  = \Psi_t E_{t}.   \end{equation}
The construction is as follows. 
Let $ t $ be any $\lambda$-tableau 
and let $ k $ be an integer between $ 1 $ and $ n $. 
The radial length between the nodes $ t[k] $ and $ t[k-1] $ is defined as 
$ h_{t,k} =h_k := r_t(k-1) - r_t(k)   $. 
Let 
$ d(t) = \sigma_{i_1} \sigma_{i_{2}}\ldots  \sigma_{i_N}     $
be a reduced expression of $ d(t) $. We associate with it 
a sequence of tableaux $ t_1= t^{\lambda}, t_2,   \ldots, t_{N+1}= t  $ by setting 
recursively $ t_{ k+1} := t_{k} s_{ i_k}  $. 
Then $ \Phi_{t} $ and
$ \Psi_{t} $ are given by the formulas 
\begin{equation}
\begin{array}{c}
\Phi_{t} :=  
\left(\sigma_{i_1} -\frac{1}{h_{t_1, i_1}} \right) 
\left(\sigma_{i_{2}} -
\frac{1}{h_{t_{2}, i_{2}}} \right) 
\ldots 
\left(\sigma_{i_N} -\frac{1} {h_{t_N, i_N}} \right)
\\
\Psi_{t} :=  
\left(\sigma_{i_1} +\frac{1}{h_{t_1, i_1}} \right) 
\left(\sigma_{i_{2}} +
\frac{1}{h_{t_{2}, i_{2}}} \right) 
\ldots 
\left(\sigma_{i_N} +\frac{1} {h_{t_N, i_N}} \right)
\end{array}
\end{equation}
As noted in [Mu92], 
$ \Phi_t $ and $  \Psi_t $ actually do depend on the chosen
decomposition of $ d(t) $, and not just on $ d(t)$, and so the notation is slightly
misleading. On the other hand, the key property (\ref{key-property}) holds independently
of the choice of reduced expression of $ d(t) $, and so we just take anyone.

\medskip
Our aim is to construct similar 
elements for $ E_T $ and $ E_S$. For this we need the following commutation rules between  
$ \sigma_k $ and the powers $L_k^m $ and $ (L_{k} -a)^m $. 
\begin{lemma}{\label{following-commutation-rules}} 
For $ m \in \mathbb N $ and $ a \in R $ the following formulas hold:
$$ \begin{array}{l}
a) \, \, \sigma_k L_{k}^m = L_{k-1}^m  \sigma_k + 
\sum_{i=0}^{m-1} L_{k-1}^i L_{k}^{ m-i-1}  \\
b) \, \, \sigma_k (L_{k} -a)^m = (L_{k-1} -a)^m  \sigma_k + 
\sum_{i=0}^{m-1} (L_{k-1} -a)^i (L_{k} -a )^{ m-i-1}.  \\
\end{array}
$$
\end{lemma}
\begin{pf*}{Proof}                                  
Formula $a) $ is proved using a straightforward induction on the 
commutation rules given in (\ref{Jucys-Murphy}). Formula $b) $ is proved 
the same way, since $ L_k -a $ satisfies the same commutation rules with $ \sigma_k $ as $ L_k$ does.
\end{pf*}                                  
We generalize the concept of radial length to tableaux classes by setting
$$ h_{T,k} = h_k := r_T(k-1) - r_T(k) \in \Z \subseteq R    $$
for $ k $ any integer between $ 1 $ and $ n $. It depends on the choices
of $ r_T(k ) $ and is therefore only unique modulo $ p $.

\medskip
We are now going to construct certain elements 
$ \Psi_{L, \,t}  $, verifying a generalization of 
(\ref{key-property}) for the $ E_T$'s.
Set first
$ h_L ( k) = h_L  := L_{k-1} - L_{k}$.
Modelled on $ \Psi_t $, we shall construct $ \Psi_{L, \,t} $ as products of expressions of the form 
$$ \sigma_k - \frac{1}{h_L}.$$
On the other hand, 
for such 
expressions to make sense in general, one would 
need to consider an appropriate completion of the group ring, and define $ \frac{1}{h_L} $ 
inside it as 
a power series. 
We here take a simpler approach, 
always considering 
$ L_{k} $ and $ L_{k-1} $ as elements of $ \End_{   \F}(V) $
for $ V $ a $\F$-vector space such that $ L_{k-1} -  L_{k} + \alpha 
\in  \End_{ { \F}} (V) $ is nilpotent 
for some $ \alpha \in \F^{\times} $.
Under that assumption, 
$ \frac{1}{h_L} $ can be defined as the corresponding geometric
series, which is finite. The next Lemma should be seen in this light.

\begin{lemma}{\label{power_series}} 
Suppose that $ s, t $ are standard tableaux with $ s= t \sigma_k  $ and that $ T:= [t] $ and $S := [s]$ are different 
tableaux classes.
Let $ h := h_{T,k} $. 
Then 
$ L_{k-1} - L_{k}  - h $ acts 
nilpotently in $ E_T (  \F S_n) $.
Especially, 
$ L_{k-1} - L_{k}$ is invertible as an element of 
$ \End_{\F}(E_T ( \F S_n)) $.
\end{lemma}
\begin{pf*}{Proof}                                  
Notice that since $ E_T \in \A $, we have that the product $ E_T (  \F S_n) $ is well defined. 
Consider first $  L_{k-1} - L_{k} - h $ as an element of $ \A $.
Using formula 
({\ref{Mu_page_506}}) we have that  
$$ ( L_{k-1} - L_{k} - h)^N = 
\sum_{ u  } \, (r_u (k-1) - r_u (k) - r_T (k-1) + r_T (k) )^N  E_u. $$ 
Multiplied by $ E_T$ it gives the formula 
$$ ( L_{k-1} - L_{k} - h)^N E_T= 
\sum_{ u \in T } \, (r_u (k-1) - r_u (k) - r_T (k-1) + r_T (k) )^N  E_u. $$ 
Each coefficient of $ E_u $ is here a multiple of $ p $.
Hence we may take
$ N $ large enough
for $ ( L_{k-1} - L_{k} - h)^N E_T$ to belong to $p^m \A $. We reduce modulo $ p $ and get 
the statement of the Lemma.
\end{pf*}

We can now prove 
the following Lemma.
\begin{lemma}{\label{the_following_lemma}} 
Let $ s,t, S, T , h, m $ be as in the previous Lemma and 
let $ h_L := L_{k-1} - L_{k} $. View $ 1/h_L $ as an element of 
$ \Hom_{ \F}(E_T ( \F  S_n),  \F S_n ) $ via 
the previous Lemma. Then
for $ N \in \mathbb N $ and $ a \in \F  $ the following formulas hold in 
$ Hom_{ \F } (E_T ( \F S_n),  \F S_n ) $
$$ \begin{array}{l}
a) \, \, (\sigma_k - \frac{1}{h_L}) L_{k}^N = L_{k-1}^N ( \sigma_k - \frac{1}{h_L}) \\
b) \, \, (\sigma_k - \frac{1}{h_L}) (L_{k}-a)^N = (L_{k-1}-a)^N ( \sigma_k - \frac{1}{h_L}).
\end{array}
$$
\end{lemma}
\begin{pf*}{Proof}                                  
Using the previous Lemma and the fact that $ E_T $ commutes with $ L_k^r $ and $ (L_k-a)^r $ 
we 
first notice that 
the expressions are well defined transformations of $ E_T ( \F  S_n) $.	
Let us now show $a)$. Since $ L_k $ and $ L_{k-1} $ commute it is equivalent to 
$$ h_L(\sigma_k L_k^N -  L_{k-1}^N  \sigma_k) = L_k^N - L_{k-1}^{N} $$
and hence, using Lemma {\ref{following-commutation-rules}}, to the valid expression
$$ (L_k -L_{k-1} ) \sum_{i=0}^{N-1} L_{k-1}^i L_{k}^{ N-i-1}  = L_k^N - L_{k-1}^{N}. $$  
Formula $ b) $ is proved the same way.
\end{pf*}                                  
We now obtain the following result.
\begin{lemma}{\label{the_following_result}} 
Let the notation be as above.  
Then we have 
$$ 
\left(\sigma_k - \frac{1}{h_L(k)} \right) E_T = E_S \left(\sigma_k - \frac{1}{h_L(k)}\right) 
$$
in 
$  \Hom_{\F} (E_T (\F S_n),  \F S_n) $.
\end{lemma}
\begin{pf*}{Proof}
The proof is obtained by combining Lemma \ref{future_use},
\ref{due_to_Murphy} and \ref{the_following_lemma}.
\end{pf*}                                  
The Lemma is a generalization of Lemma 6.2 from [Mu92], where
$ L_k, L_{k-1} $ and hence $ h_L $ act semisimply.
Note that the second minus sign 
is there a plus sign, corresponding to the fact that 
the eigenvalues of $ h_L $ on $ E_s $ and $ E_t $ 
are equal but with opposite signs.

\medskip
Set $ T^{\lambda} := [t^{\lambda}] $. 
For $ d(t) = \sigma_{i_1}  \sigma_{i_{2}}  \ldots \sigma_{i_N} $ in reduced form we define 
\begin{equation}{\label{in_reduced_form}}
\begin{array}{c}
\Psi_{L,d(t)} :=  
\left(\sigma_{i_1} -\frac{1}{h_L (i_1)} \right)
 \left(\sigma_{i_2} -\frac{1}{h_L(i_2)} \right) \ldots
 \left(\sigma_{i_N} -\frac{1}{h_L(i_N)} \right) 
\end{array}
\end{equation}
where $ \frac{1}{h_L (i_j	)}  $ is set to $ 1 $ when 
$  [ t_j  ] = [ t_{j-1}] $.
Combining Lemma {\ref{commute}} and \ref{the_following_result} we get 
the following Theorem.
\begin{Thm}{\label{we_get_that_Psi}}
$ E_{T^{\lambda}} \Psi_{L,d(t)}  =  \Psi_{L,d(t)} E_{ T}   .$
\end{Thm}

We view $ \sigma_k - \frac{1}{h_L(k)}   $ as an analogue of the Khovanov-Lauda generator
$ \psi_i $, or more precisely of the element denoted $ \phi_i $ in [BK]. 
These {\it intertwining} elements
are the starting point of their work. In our approach 
the $ \phi $-elements have a representation theoretical interpretation coming 
from the theory of the seminormal basis
whereas they appear somewhat pulled out of the sleeve in [BK].

\section{A cellular basis.}
In this section we use the results from the previous sections to construct 
a cellular basis for $ \F S_n $. 
Our construction is inspired by the one given by J. Hu and A. Mathas in [HuMa].

\medskip
Let us first introduce some
notation. 
For $ \lambda $ a partition of $ n $ we 
let $ S_{\lambda} $ denote the row stabilizer of $ t^{\lambda} $. Let $ x_{\lambda} $ and 
$ y_{\lambda} $ be the elements of $ \A	 $ given by 
$$ x_{\lambda} = \sum_{ \sigma \in S_{\lambda} }   \sigma \,\,\,\, \mbox{and} \, \, 
y_{\lambda} = \sum_{ \sigma \in S_{\lambda} } (-1)^{| \sigma | } \sigma $$
where $ | \sigma | $ is the sign of $ \sigma $. 
For a pair $ (s, t) $ of $ \lambda$-tableaux we define 
$$ x_{s t } = d(s)^{-1} x_{\lambda } d(t)  \,\,\,\, \mbox{and} \, \, 
y_{s t } = d(s)^{-1} y_{\lambda } d(t). $$
If $ s $ is a $\lambda$-tableau we get 
that $ x_{ss} $ is the sum of the elements of the row-stabilizer of $s$.	
A similar comment applies to $ y_{s s } $.

\medskip
The set $ \{ x_{s t } \} $ with $ (s, t ) $ running over pairs of standard $ \lambda$-tableaux and
$ \lambda $ over partitions of $ n $ gives 
Murphy's standard basis for $ \A$. Similarly 
$ \{ y_{s t } \} $ gives 
the dual standard basis.
They are cellular bases in the sense of Graham and Lehrer, [GL], with respect to the dominance order. This implies 
that $ \Aupe $ and $ \Aup$, defined by  
$$
\begin{array}{c}
\Aupe := \Span_{ R} \{ x_{st} | \,
(s,t) \mbox{ pair of }  \mu \mbox{-tableaux with }   
\mu \unrhd \lambda \} \\
\Aup := \Span_{ R} \{ x_{st} | \,
(s,t) \mbox{ pair of }  \mu \mbox{-tableaux with }   
\mu \rhd \lambda \}
\end{array}  
$$
are ideals of $ \A $. 
The associated left cell module is 
$$ C(\lambda) :=  R \{  x_{s \lambda}\, | \, s \mbox{ is a} \, 
\lambda \mbox{-tableau} \}  \mbox{ mod } \Aup. $$
In modern terminology it is often referred to as the Specht module, 
although it rather corresponds to the dual 
Specht module defined via Young symmetrizers. 

\medskip
Working over the ground fields $ \F $ and $ \Q$, we get ideals 
$ \AFup, \AFupe  $ and $ \AQup, \AQupe  $ of $ \AF $ and $ \AQ$, using constructions similar to the ones for 
$ \Aup, \Aupe$.
Similarly, we get cell modules $\overline{ { C(\lambda) }} $ and 
$  C_{\Q}(\lambda) $ for $ \AF $ and $ \AQ$. 
We use the same notation $ x_{s \lambda} $ for the classes of $ x_{s \lambda} $ in $ C(\lambda) $, 
$ \overline{ C(\lambda) }$ or $ { C_{\Q}(\lambda) }$. 
They form bases 
for $ C(\lambda) $ over $ R $, for $ \overline{ C(\lambda) }$ over $ \F $ and 
for $ C_{\Q} (\lambda) $ over $ \Q $, when $ s \in \Std(\lambda)$.  Hence we have the base change properties
$\overline{  C(\lambda) } := C(\lambda) \otimes_R \, \F $ and  
$  C_{\Q}(\lambda)   := C(\lambda) \otimes_R \, \Q $. 


\medskip
We need to recall another basis for $ \A $ that was also introduced by Murphy, see
[Mu92].  For $ \lambda \in \Par_n $, we define 
\begin{equation}{\label{following_murphy}}
 \xi_{\lambda } = 
\prod_{i=1}^n (L_i + \rho_{\lambda}(i) ) 
\end{equation}
where $ \rho_{\lambda}(i) = k $ for $ t^{\lambda}(i)= [k,l] $, that is 
$ \rho_{\lambda}(i) $ is the row number of the $ i $-node of $ t^{\lambda} $.
For any pair $ (s, t) $ of $ \lambda$-tableau we set 
$$ \xi_{s t} := d(s)^{-1} \xi_{\lambda } d(t). $$
Then $ \{ \xi_{s t} \} $ is a basis for $ \A $ when $ (s, t) $ runs over
the same parameter set as above. This follows from 
Theorem 4.5 of [Mu92], saying that  
\begin{equation}{\label{xi-basis}} 
\xi_{\lambda }   = x_{\lambda } +\sum_{ t \in \RStd(n), \,  t \unrhd t^{\lambda}  } \,  x_{t t }
\end{equation}
and Theorem 4.18 of [Mu95], saying that for any two tableaux $ u,v \in \RStd(n) $ of the same shape, the
element $  x_{uv} $ can be expanded in terms of standard basis elements $  x_{u^{\prime} v^{\prime}} $
where $ (u^{\prime}, v^{\prime}) \unrhd (u,v)$. We also get from this and (\ref{xi-basis}}), 
using that $ \Aup$ is an ideal in $ \A$, 
that the images of $ \{ \xi_{s \lambda} \}  $ 
in $ C(\lambda ) $
coincide with $ \{ x_{s \lambda} \}  $, for $ s $ standard $ \lambda$-tableaux.

\medskip  
Motivated by the construction done by Hu and Mathas in [HuMa]
we now introduce for each pair of standard tableaux $ (s, t ) $ of the same shape $ \lambda $
the following elements of $ \AF   $
\begin{equation}{\label{psi-elements}} \psi_{s t} := \Psi_{L, d(s)}^{\ast}
\xi_{\lambda} E_{T^{\lambda}} \Psi_{L, d(t)}
\end{equation}
where $ \Psi_{L, d(s)}, \Psi_{L, d(t)} $ are as in 
(\ref{in_reduced_form}) and where once again $ \ast$ is the antiautomorphism that 
fixes the transpositions. Since $ \frac{1}{h_L(k)} $ is a polynomial 
expression of Jucys-Murphy elements, $ \Psi_{L, d(s)}^{\ast}$ is obtained from 
$ \Psi_{L, d(s)}$ by reversing the factors.
Note that in $ \psi_{s t} $ 
the two middle factors $ \xi_{\lambda} $ and $ E_{T^{\lambda}} $ commute. 

\medskip	
We aim at proving that the set of $ \psi_{s t} $ 
is a cellular basis for $ \AF $ when
$ (s, t) $ runs over pairs of standard tableaux of the same shape.
We begin with the following preparatory Lemma.
\begin{lemma}{\label{preparatory-lemma}} 
For every partition $ \lambda $ of $ n $ we have 
triangular expansions 
$$
\xi_{\lambda} = x_{\lambda }  + \sum_{(u,v) \vartriangleright (s,t)} c_{uv} x_{uv}, \mbox{ }  \, \, \,
\xi_{\lambda} E_{T^{\lambda}} = x_{\lambda }  + \sum_{(u,v) \vartriangleright (s,t)} d_{uv} x_{uv} 
$$
where $ c_{uv} \in R,  d_{uv} \in \F$.
\end{lemma}
\begin{pf*}{Proof}
The first expansion follows from  (\ref{xi-basis}) 
and Theorem 4.18 of [Mu95].

In order to obtain the second expansion, we first recall 
Corollary 2.15 of [HuMa1] saying that for any pair of 
standard tableaux $ (s,t) $ of the same shape we have 
\begin{equation}{\label{Cor_2.15}}
x_{st} L_k = r_t(k) x_{st} + \sum_{(u,v) \vartriangleright (s,t)} c_{uv} x_{uv} 
\end{equation}
where $ (u,v) $ runs over pairs of standard tableaux of the same shape and 
$ c_{uv} \in R $.	
Especially, we get 
$$   x_{\lambda} L_k = r_{\lambda}(k) x_{\lambda} + \sum_{(u,v)  \rhd (t^{\lambda}, t^{\lambda})} 
c_{uv} x_{uv}. $$
From this we get via the definition of $ E_t $ that 
$$   x_{\lambda} E_{t } =  \delta_{ t t^{\lambda}} x_{\lambda} + \sum_{(u,v)  \rhd (t^{\lambda}, t^{\lambda})}
d_{uv} x_{uv} $$
for certain $ d_{uv} \in \Q$ where $ \delta_{ t t^{\lambda}} $ is the Kronecker delta. 
Using the definition of $ E_T $ we deduce from this that 
$$   x_{\lambda} E_{T } =   x_{\lambda} + \sum_{(u,v)  \rhd (t^{\lambda}, t^{\lambda})}
d_{uv} x_{uv} $$
where $ (u,v) $ still is a pair of standard tableaux of the same shape but 
where we may now assume that $ d_{uv}  \in R$. We then finish the proof of the Lemma by reducing modulo $ p$.
\end{pf*}
{\bf Remark}.
For $ \mu := Shape(\sigma) $ and $ \nu := Shape(s)  $ we define $ (\sigma, \tau)  \succeq (s, t) $ 
if either $ (\sigma, \tau)  \unrhd (s, t)  $ and $ \mu = \nu $, or if $ \mu \rhd \nu$. Then 
the triangularity property ({\ref{Cor_2.15}}) appears already in [AMR] if $ \unrhd $ is replaced 
by $ \succeq$.

\medskip
The next result gives the promised cellularity property for $ \{  \psi_{s t} \} $. 
\begin{Thm}{\label{cel-basis}} 
For pairs of standard tableaux $ (s, t) $ of the same shape we have 
$$ \psi_{st}  = x_{st} + 
\sum_{ (\sigma, \tau)  \succ (s, t)}  
a_{ \sigma \tau} x_{ \sigma \tau} , \, \, \, 
a_{ \sigma \tau} \in  \F.
$$ 
where $ \succeq  $ is the partial order introduced in the above remark.
Moreover, with respect to the dominance order, the set 
$$ \{  \psi_{s, t} \, | (s, t)   \mbox{ pair of standard tableaux of the same shape} \} $$
defines a cellular basis for $ \AF  $ with cell modules $ C(\lambda) $.
\end{Thm}
\begin{pf*}{Proof}
We have
$$ \psi_{s t} := \Psi_{L, d(s)}^{\ast} E_{T^{\lambda}} 
\xi_{\lambda}   \Psi_{L, d(t)}. $$
For $ d(s) = \sigma_{i_1}  \sigma_{i_{2}} \ldots \sigma_{i_N} $
in reduced form, 
the rightmost term of $ \Psi_{L, d(s)}^{\ast} $ is 
$ \sigma_{i_1} - \frac{1}{h_L(i_1) } $. 
Let us consider 
its action on $ E_{T^{\lambda}} 
\xi_{\lambda} $. Assume first that $ \frac{1}{h_L(i_1) } \neq 1 $. 
Now $ \frac{1}{h_L(i_1) } = \frac{1}{L_{i_1 -1} -L_{i_1}} $ 
is a linear combination of terms of the form
$$ (L_{i_1 -1} -L_{i_1} -r_{{\lambda}}(i_1 -1) + r_{{\lambda}} (i_1) )^l $$
and hence by ({\ref{Cor_2.15}}) it acts upper triangularily. Thus $ \sigma_{i_1} - \frac{1}{h_L(i_1) } $
maps $ E_{T^{\lambda}} \xi_{\lambda} $ to $ x_{ \sigma_{i_1} t^{\lambda} , \lambda } $
plus higher terms with respect to $ \succeq$. If $ \frac{1}{h_L(i_1) } = 1 $ the same conclusion holds trivially.
We repeat this argument 
for the other terms of $ \Psi_{L, d(s)}^{\ast} $ and then for the terms of 
$ \Psi_{L, d(t)}$, and obtain the triangularity statement of the Theorem.

\medskip
From this we deduce that 
$$  \AFup = \Span_{ \F} \{ \, \psi_{st} \, | \,
(s,t) \mbox{ pair of }  \mu \mbox{-tableaux with }   
\mu \rhd \lambda \}. $$
An argument similar to the one given in Theorem 5.8 of [HuMa] 
now gives 
the cellularity of $ \{ \psi_{st} \} $ 
with $\ast$-involution satisfying $ \psi_{st}^{\ast} = \psi_{ts}$.
\end{pf*}

From the general theory of cellular algebras there is an associated bilinear invariant 
form on $ \overline{C(\lambda )} $, that 
we denote
$ \langle \cdot, \cdot \rangle_{\lambda} $. It it given by 
$$ \psi_{ \lambda s} \psi_{t \lambda} = 
 \langle \psi_{s \lambda}, \psi_{t \lambda} \rangle_{\lambda} \, \psi_{\lambda} \mbox{ mod } 
\AFup$$ 
Its radical $ \rad_{\lambda} $ is a submodule of $ \overline{C(\lambda )} $ and 
$ \overline{ C(\lambda )}/\rad_{\lambda} $ is either simple or zero. 
An important point of the theory of cellular algebras is that this gives rise to a classification of the 
simple modules for $ \AF $. Indeed, every simple module for $ \AF $ is of the form
$ \overline{ C(\lambda )}/\rad_{\lambda} $ for a unique $ \lambda$. 

\medskip
Our next Lemma shows that 
$ \langle \cdot, \cdot \rangle_{\lambda} $ is in block form with respect to our basis.
Note that bases that block diagonalise the bilinear form have also been found in [MaSo] and in [BKW]. 
\begin{lemma}{\label{diagonalizing}} 
The basis $ \{ \psi_{s \lambda} | \, s \mbox{ standard } \lambda \mbox{-tableau} \} $ 
of $ \overline{C(\lambda )} $ is in block form with 
respect to $ \langle \cdot, \cdot \rangle_{\lambda} $
with blocks given by the tableaux classes.
\end{lemma}
\begin{pf*}{Proof}
Suppose that $ s, t $ are standard $ \lambda $-tableaux and that the tableau classes 
$ S:=[s] $ and $ T:=[t ] $
are different.
Then we have that 
$$ 
\psi_{ \lambda s} \psi_{t \lambda} =  
 \xi_{\lambda}  E_{T^{\lambda}}
\Psi_{L, d( s) } \Psi_{ L ,d(t) }^{\ast} E_{T^{\lambda}}   \xi_{\lambda}.
$$
Using Theorem {\ref{we_get_that_Psi}} and its $ \ast$-version, and noting that $ E_{U }^{\ast} = 
E_{U} $ for all $ U $ since $ E_{U}$ is a sum of products of Jucys-Murphy operators, we get that 
this is equal to
$$ 
\xi_{\lambda} \Psi_{L, d( s) }  E_{S}
 E_{T}  \Psi_{ L d(t) }^{\ast}  \xi_{\lambda}.
$$
But $ E_{S}
 E_{T} = 0 $ and the Lemma follows.
\end{pf*}

\section{Specht modules and Jucys-Murphy operators.}
In this section
we give a new realization of the Specht modules, 
using 
Jucys-Murphy operators. 

\medskip
An essential ingredient of our construction 
is the use of what we denote the Gelfand-Zetlin subalgebra of $ \A $
as a kind of Cartan subalgebra of a semisimple Lie algebra. This is 
in accordance with ideas  
promoted by Okounkov and Vershik in the article ``A new approach to 
the representation theory of the symmetric group'', [OV]. 
Their approach also applies to a wider class of algebras than the group algebra
of the symmetric group, but relies heavily on the algebras being semisimple.
Moreover, the Specht modules themselves 
have in their approach apparently so far only been treated from the ``old'' point of view. 
In this section we realize the Specht module as induced modules from the Gelfand-Zetlin 
subalgebra, at least over $ R $ and $ \Q$ and thus partially remedy these deficiencies.
It would be interesting to investigate to what extent 
these results hold in positive characteristic.

\medskip
Define $ \GZ  \subseteq   {\A}  $ to be 
the Gelfand-Zetlin algebra, the $R$-subalgebra of $ \A $ generated by 
the Jucys-Murphy operators: $$ \GZ := \langle L_i \, | \, i= 1, \ldots , n \rangle. $$
This definition is not quite equivalent to the one used by for example Okounkov and Vershik in 
[OV]. They first of all work over a field of characteristic zero and even in that case, 
our definition of the Gelfand-Zetlin algebra is actually a Theorem in [OV] that characterizes 
the subalgebra.

\medskip
$ \GZ $ is a commutative subalgebra of $ {\A} $ and it
contains the center $ Z(\A)  $ of $ {\A} $ -- indeed by Theorem 1.9 of [Mu83] we know that
$ Z(\A)  $ 
consists of the symmetric polynomials in the $ L_k$.

\medskip
We aim at defining an induction functor from $ \GZ$-modules to $ \A$-modules. 
For this we first need to state 
a few categorical generalities on $ R$-modules.

\medskip
For an $ R$-module $ M $ we define $ M^{\ast} := \Hom_{R}(M, R)$. 
If $ M $ is also a left $ \A $-module,
$ M^{\ast} $ becomes a right $ \A $-module and vice-versa.
Let $R \modfg$ denote the category of 
finitely generated $ R$-modules and 
let $\A \modfg$ denote the subcategory whose objects are also left $ \A$-modules.

\medskip
Since $ R $ is Euclidean, we have for $ M \in R \modfg $ that $$ M = \Fr(M) \oplus \T(M) $$
where $ F(M) $ is the free part of $ M $ and $ \T(M) $ the torsion part of $M$. 
If $ f: M \rightarrow N $ is a morphism in $R \modfg$ then clearly $ f(\T(M)) \subseteq \T(N)$ and from this we deduce 
that $ M \mapsto \T(M) $ 
is a left exact functor on $R \modfg$. On the other hand 
$ M \mapsto \Fr(M) $ is an exact functor. Indeed, we may define it as 
$\Fr(M) := M^{\ast \ast} $ which shows that it is a covariant functor in the first place. 
But for $ M \in R \modfg $ the canonical map $ M \rightarrow M^{\ast \ast} $
induces an isomorphism $ M/\T(M) \rightarrow \Fr(M)  $. This gives a natural transformation 
from the functor $ M \mapsto M/\T(M) $ to $ \Fr$ and hence, since $M \mapsto  M/\T(M)   $ is 
right exact, we get that $ \Fr $ is right exact as well, whereas left exactness follows 
directly from the definitions. 

\medskip
From this we get that $ \Fr $ induces an exact functor on $\A \modfg$. 
Indeed, if $ M$ is a left $ \A $-module then
also $ \Fr(M) = M^{\ast \ast } $ is a left $ \A $-module and exactness follows from exactness at 
$R \modfg$ level.

\medskip
We let $\GZ \modfg$ denote the subcategory of $R \modfg$ whose objects are also 
$ \GZ$-modules. Finally, we define 
$R \modf$ as the category of 
finitely generated free $ R$-modules and 
$\A \modf$ as the subcategory whose objects are also left $ \A$-modules.

\medskip
After these preparations we are in position to define the induction functor.
For $ M \in \GZ \modfg $ we define 
$$ \Ind(M) := \Fr ({  \A} \otimes_{\GZ} M ).$$ 
Then  $\Fr ({  \A} \otimes_{\GZ} M )
\in \A \modf $.
Furthermore, by the above considerations we have that $ M \mapsto \Ind(M) $ 
is a right exact functor from $\GZ \modfg$ to $\A \modf$.

\medskip
An important property of $ \Ind $ 
is the following Frobenius reciprocity rule
$$ \begin{array}{c}
\Hom_{\GZ}(M, N ) \cong  \Hom_{{\A}}(\Ind(M), N  )
\end{array}
$$ 
for $ M \in \GZ \modfg $ and $ N \in \GZ \modf $. It follows from 
$$ \Hom_R(M, N ) \cong \Hom_R(F(M), N )  $$
for $ M \in R \modfg, \, N \in R \modf  $
and the usual 
Frobenius reciprocity for induction.

\medskip
For us the most important case of the above construction is the following. 
Let $ \lambda $ be a partition of $ n $ and
let $ I_{\lambda} $ be the ideal of $ \GZ $ generated by $ L_i -r_{\lambda}(  i)$
for $ i = 1, \ldots, n $. 
Set $$ 1_{\lambda} := \GZ / I_{\lambda}. $$ 
Then we may consider $ 1_{\lambda} $ as a left $ \GZ $-module. 
As we point out in the final remarks of this section, 
it is free of rank one over $ R $ with generator $ 1 $. The action 
of $ L_i $ on $ 1$ is multiplication by $ r_{\lambda}(  i)$. 
We next define 
$$ \Ind(\lambda) := \Ind(1_{\lambda}). $$
We aim at studying $  \Ind(\lambda)$ at some depth, 
our main result being a proof of the isomorphism
$ \Ind(\lambda) \cong  C(\lambda) $.
The 
following Lemma is a first step towards this.

\medskip
Let $ t_{\lambda} $ be the lowest $ \lambda$-tableau having $ 1,2 \ldots, n $ 
filled in along columns and define $ s_{\lambda}:= d(t_{\lambda} ) \in  S_n $. Set 
\begin{equation}{\label{z-element}}
z_{\lambda} := x_{\lambda} s_{\lambda} y_{\lambda^{\prime}} =
x_{\lambda} s^{ \,-1}_{\lambda^{\prime}} \, y_{\lambda^{\prime}}
 \in \A.
\end{equation}
Then $ x_{\lambda} \A $ is isomorphic to the right permutation module studied in [J] where the 
isomorphism maps the tabloid $ \{ t^{\lambda} \} $ to $ x_{ \lambda} $. The alternate column 
sum of $ t^{\lambda} $ is $ \kappa_{ t^{\lambda}} = 
s^{ \,-1}_{\lambda^{\prime}}  \, y_{\lambda^{\prime}} s_{\lambda^{\prime}} = 
s^{}_{\lambda} \, y_{\lambda^{\prime}} s_{\lambda}^{-1} 
$ 
and so the isomorphism maps the polytableau $ e_{t^{\lambda}} $ of [J] to 
$ x_{\lambda} s^{ \,-1}_{\lambda^{\prime}} \, y_{\lambda^{\prime}} s_{\lambda^{\prime}}$.
In other words $ x_{\lambda} s_{\lambda} y_{\lambda^{\prime}} \A $ identifies with the Specht module 
$ S^{\lambda} $ of [J].
For $ s $ a $ \lambda$-tableau we set 
$$ z_{\lambda s} := x_{\lambda} s_{\lambda} y_{\lambda^{\prime}} d(s^{\prime}). $$
Then $ \{ z_{\lambda s} \, | \, s \in \Std(\lambda) \} $ is a basis for $ z_{\lambda} \A $, 
see eg. Lemma 5.1 and Theorem 5.6 of [DJ].
\begin{lemma}{\label{with_the_above_notation}}
With the above notation we have 
$$   \{ x \in {\A} \, | \, L_i x = r_{\lambda}(i) \, x  \mbox{ for all } i \}  =  z_{\lambda} \A.  $$
\end{lemma}
\begin{pf*}{Proof}
Let us denote by $ { }_{\lambda}{\cal LS}_n $ the left hand side of the Lemma. We first prove 
that $  { }_{\lambda}{\cal LS}_n   \supseteq z_{\lambda} \A $. Now 
$ { }_{\lambda}{\cal LS}_n $ certainly is a right 
submodule of $ \A $ and 
it is known, see for example [Mu92] page 498, that 
\begin{equation}{\label{we_know_from_[Mu92]}} x_{ab } s_{\lambda} y_{\lambda^{\prime}} = 0 \, \, \, 
\mbox{ unless } \mu \unlhd \lambda 
\end{equation}
where $ \mu = Shape(s) = Shape(b) $.
Combining this with the fact that $ L_i $ acts upper triangularily on the 
$ \{ x_{st} \} $-basis, as is seen by applying $ \ast $ to (\ref{Cor_2.15}), we 
find that $ z_{\lambda}  $ 
belongs to $ { }_{\lambda}{\cal LS}_n $, from which the inclusion $ \supseteq $ indeed follows. 


\medskip 
In order to show the other inclusion $ \subseteq $ we first work in $ \AQ$, and define
$$ { }_{\lambda}{\cal LS}_{\Q,n} : =  \{ x \in \AQ \, | \, L_i x = r_{\lambda}(i) \, x  \mbox{ for all } i \}. $$
Setting $ t = t_{\lambda}  $ we recall from [Mu92] page 511 that 
\begin{equation}{\label{page_511}}
E_{\lambda} = h_{\lambda}^{-1}  z_{\lambda t} \Psi^{\ast}_t 
\end{equation}
and from this we deduce that $ { }_{\lambda}{\cal LS}_{\Q,n}  = z_{\lambda} \AQ $. 
We then get that 
$$ { }_{\lambda}{\cal LS}_n = 
z_{\lambda} \AQ
\cap \A = z_{\lambda} \AQ \cap x_{\lambda} \A.  $$
Here the last equality follows from the facts that $ \{ x_{st} \}$ is an $ R$-basis 
of $ \A $ and that $  z_{\lambda} \AQ \subseteq x_{\lambda} \AQ $.
Finally, since $  z_{\lambda} \AQ  = S_{\Q}(\lambda)$ is the Specht module defined over $ \Q$, 
we get from Corollary 8.9 of [J], which is based on the Garnir relations, 
that 
$$ z_{\lambda} \AQ \cap x_{\lambda}\A \subseteq z_{\lambda} \A $$
and the Lemma is proved.
\end{pf*}

Recall that $ \A  $ is equipped with a 
symmetric nondegenerate bilinear form $ \langle \cdot, \cdot  \rangle $, given by 
$$ \langle a, b \rangle := \coeff_1 (ab) $$ 
where $ \coeff_1 (x) $ is the coefficient of $ 1$ when $ x \in \A $ is expanded in the 
canonical basis of $ \A $. 
It is associative in the following sense
$$ \langle ab, c  \rangle = \langle a, bc  \rangle \,\, \,  
 \mbox{ for all } a,b,c \in  \A. $$
The form induces an $ \A$-bimodule isomorphism 
$ {\A}  \cong \Ad = \Hom_{R}({\A}, R) $ where the $ \A$-bimodule structure on 
$ \Ad$ is given as follows
$$ a f(x) b := f(bxa) \, \, \, \, \, \, 
 \mbox{ for all } a,b, x \in  {\A}   \mbox{ and }  f \in  \Ad .
$$
We can now prove the promised result on $ \Ind(\lambda)$. 
\begin{Thm}{\label{promised}}
For $ \lambda $ any partition of $ n $ there is an isomorphism of $ \A $-modules
$$   \Ind(\lambda) \cong C(\lambda). $$
\end{Thm}
\begin{pf*}{Proof}
Define $ \LGZ_{\lambda } := \sum_i {\A } (L_i -r_{\lambda}( i )) $.  Then 
$ \LGZ_{\lambda } $ is a left ideal of $ \A $ and by the definitions we have that 
\begin{equation}{\label{and_by_the_defin}}
 \Ind(\lambda) = \Fr({\A } / \LGZ_{\lambda }) = ({\A } / \LGZ_{\lambda })^{\ast \ast}. 
\end{equation}
But $ \langle \cdot , \cdot \rangle $ is nondegenerate, and therefore it induces an isomorphism of 
the right $ \A $-modules
$$ ({\A } / \LGZ_{\lambda })^{\ast} \cong 
(\LGZ_{\lambda })^{\perp} 
$$
where 
$ (\LGZ_{\lambda })^{\perp} := \{ x \in {\A} \, | \,  \langle x , \LGZ_{\lambda } \rangle =0 \,\}$.
On the other hand, using the symmetry, associativity and nondegeneracy of $ \langle \cdot, \cdot \rangle $ we find that 
$ x \in (\LGZ_{\lambda })^{\perp } $ iff $ (L_i -r_{\lambda} (i )) x = 0 $ for all $ i $. 
We then deduce from the previous Lemma that 
$$ (\LGZ_{\lambda })^{\perp } = z_{\lambda } \A .$$
Thus, we are reduced to showing that $ z_{\lambda } \Ad \cong C(\lambda) $. This is a little variation of 
a well-known fact, that normally is presented using either two left or two right modules.
In our setting, with one left and one right module, the pairing 
$ z_{\lambda }\A  \times   C(\lambda)     \mapsto R $ is given 
by the rule $ ( z_{\lambda t} ,  x_{s \lambda}    ) \mapsto {\coeff}_{\lambda} (z_{\lambda t}   
x_{s \lambda}) $ where for any $ u \in \A $ we define $ {\coeff}_{\lambda} (u) $ as 
the coefficient of $ x_{\lambda} $ when $ u $ is expanded in the $ x_{st} $-basis.

\end{pf*}

We now deduce the following universal property of $ C(\lambda) $. We consider it analogous to the universal property 
for the Weyl module of an algebraic group, which is a consequence of the Kempf's vanishing Theorem of 
the cohomology of the line bundle on the flag manifold given by a dominant weight, see eg. [A2], [RH3].
\begin{Thm}{\label{universal}} 
Let $ M \in \A \modf$. 
Let $$ { }_{\lambda}M := \{ m \in M \,|\, L_i m = r_{{\lambda} } (i ) m  \mbox{ for all } i \}. $$
Then 
$ \Hom_{{\A}} (  C(\lambda), M )= { }_{\lambda}M$. 
\end{Thm}
\begin{pf*}{Proof}
Any $ m \in { }_{\lambda}M $ induces a map in 
$ \Hom_{GZ_n} (  1_{\lambda}, M ) $ and then by Frobenius reciprocity a map in 
$\Hom_{{\A}} (  \Ind(\lambda), M )=  \Hom_{{\A}} (  C(\lambda), M )   $.
On the other hand, any element of $ f \in \Hom_{GZ_n} (  1_{\lambda}, M ) $ gives rise to 
an element of $ { }_{\lambda}M $, namely the image $ f(1)$.
\end{pf*}

\medskip
\noindent
{\bf Remark}.
Let $ \GZQ $ be the original Gelfand-Zetlin algebra introduced in [OV], that is the $ \Q$-subalgebra of 
$ \AQ$ generated by the Jucys-Murphy elements $ L_i, i =1, \ldots, n$. We 
let $ I_{\Q, \lambda} $ denote the ideal of $ \GZQ $ generated by the $ L_i -r_{\lambda}(  i) $ for $ i =1, \ldots, n$
and set 
$ 1_{\Q, \lambda} : = \GZQ/ I_{\Q, \lambda}$ and can then define  
$ \Ind_{\Q}(\lambda) $ as the $ \AQ $-module induced up from $ 1_{\Q, \lambda}$.
Now the same series of 
arguments as the one used above, even with some simplifications
in Lemma \ref{preparatory-lemma}, leads to the isomorphism 
$$ \Ind_{\Q}(\lambda) \cong C_{\Q}(\lambda). $$
But in this case the result could actually also have been obtained as follows.
From ({\ref{Mu_page_506}}) we have that 
\begin{equation}{\label{dim-GZ}}
 \GZQ = \langle E_t \,| \,
t \in \Std(n)  \rangle .
\end{equation}
Here the $ \{ E_t \} $ even form a $\Q $-basis for $ \GZQ$ since they are orthogonal idempotets. 
On the other hand, since $ E_{\lambda } \xi_{\lambda } = \gamma_{\lambda} E_{\lambda } $,
as is proved on page 508 of [Mu92], we get that 
the basis for $ \AQ $ constructed in the previous section, 
in this case takes the form
$$ \{ \Psi_{s}^{\ast} E_{\lambda} 
\Phi_{t} \, | \, s,t \in \Std(\lambda), \lambda \in \Par_n \}. $$
Let $ ev_{\lambda}:  \GZQ  \rightarrow 1_{\Q, \lambda} $ be the quotient map. 
Then one checks, using the fact that the contents $ r_t(i) $ determine 
$t$ uniquely, that 
\begin{equation}{\label{rankone}}
 ev_{\lambda}(E_t) = \left\{ \begin{array}{cc} 1 & \mbox{if} \, \, t = t^{\lambda} \\
0 & \, \, \, \, \mbox{otherwise.} \end{array} \right. 
\end{equation}
It now follows from $ E_{\lambda} \Phi_t =  \Psi_t E_{t} $ that 
$ \Ind_{\Q}(\lambda) $ has basis
$$ \{ \Psi_{s}^{\ast} E_{\lambda} 
 \, | \, s \in \Std(\lambda), \lambda \in \Par_n \}$$
and the claim of the Remark follows from this.

\medskip
\noindent
{\bf Remark}.
From (\ref{rankone}) we also get that $ 1_{\Q, \lambda} $ is of dimension one over $ \Q$
with basis $\{ E_{\lambda} \} $, and from this we conclude that $ 1_{ \lambda} $ is free of rank one over $R$.
Indeed,  using $  1_{ \lambda}  = \GZ/ I_{ \lambda} $ one checks that 
$  1_{ \lambda} $ is cyclic over $R$, generated by $1$.
On the other hand, using exactness of the tensor functor $ M \rightarrow M \otimes_R \Q  $ 
we get that $ 1_{ \Q, \lambda} = 1_{  \lambda} \otimes_R Q $. Hence $ \dim_{\Q} (1_{  \lambda} \otimes_R Q)  = 1 $ and 
so $  1_{ \lambda} $ is torsion free, and thus free of rank one over $R$.

\medskip
\noindent
{\bf Remark}. In general $ \A $ is {\it not} free over $ \GZ$. Indeed, if 
$ \A $ were free over $ \GZ$ then $ \AQ $ would be free over $ \GZQ$.
Since $ \{ \,E_t \,|\, t \in \Std(n) \} $
is a basis of $ \GZQ $ we can determine the dimension of $ \GZQ$.
For instance, for $ n= 3 $ we find $ \dim \GZQ = 4 $ which does not divide $ \dim \AQ = 6$.

\section{Simples.}
Let $ G $ be an algebraic group over an algebraically closed field $ k $
of characteristic $ p $.
Let $ B $
be a Borel subgroup of $G$ with maximal 
torus $ T \subset B $ and let $X(T)$ 
(resp. $  X(T)^{+}$) be the set of  weights (resp. dominant weights) 
with respect to $ B$ and $T$. For $ \lambda \in X(T)^{+} $ there is an associated
Weyl module $\Delta(\lambda) $ with unique simple quotient $ L(\lambda) $. 
It is the reduction modulo $p $ of a $ \Z $-form for a module for the corresponding complex group.
The finite 
dimensional simple modules for $ G $ are classified by $ L(\lambda) $ 
where $ \lambda \in X(T)^{+} $. We write    
$ \nabla(\lambda) := \Delta(\lambda)^{\ast}  $
where $ \ast $ is the contravariant duality functor on finite dimensional $ G$-modules. 
We may realize $ \nabla(\lambda)$ as 
the $ G $-module $ H^0(\lambda) $ of global sections of
the line bundle on $ G/B $ associated with $ \lambda $.

\medskip
Let $ \langle \cdot, \cdot \rangle_{\lambda} $ be a nonzero contravariant form on $ \Delta(\lambda) $. 
It induces a $ G $-linear map $ c_{\lambda}: \Delta(\lambda) \rightarrow 
\nabla(\lambda) $. As a matter of fact, since $ \langle \cdot, \cdot \rangle_{\lambda} $ is unique 
up to multiplication by a nonzero scalar, we have that $ c_{\lambda} $ 
generates $ \Hom_G(\Delta(\lambda), \nabla(\lambda)) $ 
and that $ im \, c_{\lambda } $ is isomorphic to $ L(\lambda) $. In this sense, $ \Delta(\lambda) $ 
and $ \nabla (\lambda) $ give rise to a {\it realization} of $ L(\lambda) $. 

\medskip
In this section we try to carry over this realization of the simple $G$-modules
to the case of the symmetric group. As we shall see, the results of the previous section provide
a suitable solution to this problem.

\medskip
Let $ M $ be a left $ \A $-module. The contragredient dual $ M^{\circledast} $ of $ M $ 
is defined to be $ M^{\ast}:=\Hom_R(M, R) $ with $ \A$-action given by $ (\sigma f)(x) := f( \sigma^{-1}x) $
for $ \sigma \in S_n,  x \in M $ and $ f \in M^{\ast} $. It is a left $ \A $-module 
as well. 

\medskip
Using Theorem 5.3 of [Mu95], with a small modification since we are working with left modules, 
we have that 
the contragredient dual of $ C(\lambda ) $ is 
\begin{equation}{\label{contra}} 
C(\lambda )^{\circledast} = 
\A \, y_{\lambda^{\prime}} s_{\lambda}^{-1} x_{\lambda}=
\A \, y_{\lambda^{\prime}} s_{\lambda^{\prime}} x_{\lambda}. 
\end{equation}
This isomorphism is also valid in the specialized situation
\begin{equation}{\label{contra}} 
\overline{C(\lambda )}^{\, \circledast} = \AF	 \, 
y_{\lambda^{\prime}} s_{\lambda^{\prime}} x_{\lambda}. 
\end{equation}
Let $ ( \cdot, \cdot)_{\lambda} $ be the bilinear form on $ C(\lambda) $ associated 
with Murphy's standard basis, following [Mu95] or the general cellular algebra 
theory, see [GL]. It is given by 
$$ (  x_{s \lambda},  x_{t \lambda})_{\lambda} = 
\coeff_{\lambda}(x_{ \lambda s } \, x_{t \lambda} )
$$
where once again $ \coeff_{\lambda}(u)  $ is the coefficient of $ x_{\lambda} $ when $u$ is
expanded in the $ x_{st} $-basis. It induces an $ \A $-homomorphism 
$ c_{\lambda}: C(\lambda) \rightarrow 
C(\lambda )^{\circledast} $, or setting 
$ z_{\lambda}^{\prime} := y_{\lambda^{\prime}} s_{\lambda^{\prime}} x_{\lambda} $ and 
using ({\ref{contra}}) and Theorem {\ref{promised}}
$$ c_{\lambda}: \Ind(\lambda) = \Fr({  \A} \otimes_{\GZ} 1_{\lambda}  )
\rightarrow \A z_{\lambda}^{\prime}.   $$
In general $c_{\lambda} $ is injective since $ ( \cdot, \cdot)_{\lambda} $ is nondegenerate 
over $ R $, but not surjective.
We can now state and prove our main result.
\begin{Thm}{\label{simples}} 
a) There is $ a_{\lambda} \in \Q $ such that $ a_{\lambda} E_{\lambda} \in \A $ and such that 
$ c_{\lambda} $ corresponds to 
$ 1 \mapsto a_{\lambda} E_{\lambda} $ under Frobenius reciprocity.
\newline
b) The simple $\AF $-module $ D(\lambda) $ associated with 
$ \lambda $  is given by $ D(\lambda) =    \AF a_{\lambda}
E_{\lambda} $.
\end{Thm}
\begin{pf*}{Proof}
By Theorem {\ref{promised}} and {\ref{universal}}, and the fact that $ C(\lambda )^{\circledast} $ 
is free over $ R $, we have 
$$ 
 \Hom_{{\A}} (  \Ind(\lambda), C(\lambda )^{\circledast})
= \Hom_{{\GZ}} (  1_{\lambda},C(\lambda )^{\circledast}  ) $$
hence $ c_{\lambda} $ is given by $ 1 \mapsto m_{\lambda} $ 
where $ m_{\lambda} \in 
{ }_{\lambda}(\A z_{\lambda}^{\prime}) $.
From Lemma 
{\ref{with_the_above_notation}} we then have 
$$ m_{\lambda} \in   \A z_{\lambda}^{\prime} \cap  z_{\lambda} \A = 
\A s_{\lambda} z_{\lambda}^{\prime}  \cap  z_{\lambda} s_{\lambda}^{-1} \A.
$$
But the Young preidempotent $ e := z_{\lambda} s_{\lambda}^{-1}  $ satisfies 
$ e^2 = \gamma_{\lambda} \gamma_{\lambda^{\prime}}  \, e $ since it can be rewritten as
$ e = x_{\lambda} \kappa_{ t^{\lambda}} $ with $ \kappa_{t^{\lambda}} $ as above.
Hence we get 
$$ m_{\lambda} = \frac{1}{\gamma_{\lambda} \gamma_{\lambda^{\prime}}} \, z_{\lambda} s_{\lambda}^{-1} m 
s_{\lambda} z_{\lambda}^{\prime} = \frac{1}{\gamma_{\lambda} \gamma_{\lambda^{\prime}}} \,
x_{\lambda} s_{\lambda} y_{\lambda^{\prime}} s_{\lambda}^{-1} m 
s_{\lambda} y_{\lambda^{\prime}} s_{\lambda}^{-1}   x_{\lambda} $$
for some $ m \in \A$. On the other hand, it is known that 
the $ R $-module $ x_{\lambda} \A \, y_{\lambda^{\prime}} $ is free of rank one, generated by 
$ x_{\lambda} s_{\lambda} y_{\lambda^{\prime}} $, see for example [Mu92] page 498, and so  
we may rewrite $ m_{\lambda} $ as follows
$$ m_{\lambda} = a_{\lambda}  
x_{\lambda} s_{\lambda} y_{\lambda^{\prime}}  s_{\lambda}^{-1}   x_{\lambda} $$
for some $ a_{\lambda} \in \Q$. We now recall the expression 
for $ z_{\lambda t} $ given on 
page 511 of loc. cit. which in our notation becomes 
$$ x_{\lambda} s_{\lambda} y_{\lambda^{\prime}} s_{\lambda}^{-1} = b_{\lambda} E_{\lambda} s_{\lambda}
E_{t} $$ where $ b_{\lambda} \in \Q $ and $ t $ is the lowest $ \lambda$-tableau. 
Applying $ \ast $ to it we get 
$$  s_{\lambda} y_{\lambda^{\prime}}  s_{\lambda}^{-1}  x_{\lambda} = b_{\lambda} 
E_{t}  s_{\lambda}^{-1} E_{\lambda}. $$
Combining these expressions and using that $ y_{\lambda^{\prime}} $ is a preidempotent, 
we find the following formula for $ m_{\lambda} $, up to a scalar in $ \Q $
$$ m_{\lambda} = E_{\lambda} s_{\lambda}
E_{t} s_{\lambda}^{-1} E_{\lambda}. $$
We then finally use the version of Young's seminormal form that is 
developed on page 152 of [RH1] and obtain 
$$ m_{\lambda} = a_{\lambda} E_{\lambda}  $$
where $ a_{\lambda} $ is a (new) scalar in $\in \Q $. This finishes the proof of $a)$.

\medskip
We next show 
$ b) $. From the definitions we have that 
$$ \A / \LGZ_{\lambda} = \A \otimes_{ \GZ} 1_{\lambda} = \Ind(\lambda) \oplus \T( \A \otimes_{ \GZ} 1_{\lambda})  $$
and so by 
$ a) $ we have that 
$ c_{\lambda}: \Ind(\lambda)   \rightarrow  C(\lambda )^{\circledast} $ 
is given by $ w \in \A  \mapsto a_{\lambda} w  E_{\lambda} $ since
$ C(\lambda )^{\circledast} $ is torsion-free.
Reducing $ c_{\lambda} $ modulo $p$ 
we get the homomorphism $ \overline{c_{\lambda}} $ 
$$  \overline{C(\lambda)} = \Ind(\lambda)  \otimes_R  \F
\stackrel{\overline{c_{\lambda}}}\longrightarrow C(\lambda )^{\circledast} \otimes_R  \F = 
\overline{\,C(\lambda)}^{\,\circledast} $$
given by 
$w \otimes 1 \mapsto a_{\lambda} w  E_{\lambda} \otimes 1$ for $ w \in \A$. 
We deduce from this that the 
image of $ \overline{c_{\lambda}} $ is the submodule of 
$ \overline{\,C(\lambda)}^{\,\circledast}$ generated by
$ \overline{a_{\lambda}E_{\lambda}} =a_{\lambda}E_{\lambda} \otimes_R \F $. 
But from the general principles explained above, this is equal to $ D(\lambda) $. The Theorem is proved.
\end{pf*}

\noindent
{\bf Remark}. So far we do not have an exact formula for $ a_{\lambda} $.
On the other hand, since
$ c_{\lambda} $ is unique up to multiplication by an element of $ R $, and since 
$ \overline{\,c_{\lambda}} $ is nonzero iff $ \lambda $ is $p$-restricted, we may simply 
choose for $ a_{\lambda} $ the least common multiple of the denominators of 
the coefficients of $ E_{\lambda} $ when expanded in the canonical basis of $ \A$.
The case where $ \lambda $ is not $p$-restricted 
is not relevant for us, of course.

\medskip
\noindent
{\bf Remark}. The Theorem gives rise to an algorithm for calculating $ \dim D(\lambda) $
that goes as follows. Let $ {\bf D}(\lambda) $ be the $ \dim S(\lambda) \times n! $ matrix over $ \F$ that has 
$ \overline{a_{\lambda} E_{\lambda}}$ in the first row 
and $  \overline{d(t)^{-1} \, a_{\lambda} E_{\lambda}}$
for $ t \in \Std(\lambda) \setminus \{t^{\lambda} \}$ in the other rows. Then $ \dim D(\lambda)= \rank {\bf D}(\lambda) $.
Note that $ E_{\lambda} $ can be calculated using formula (\ref{page_511}).
We have implemented this algorithm in the GAP system. 
We have checked $ n < 8$ for all relevant primes and found complete
match with the known dimensions for $ D(\lambda) $, as given 
by Mathas's Specht-package. 

On the other hand, although the first row of the matrix $ {\bf D}(\lambda) $ easily gives the rest of it, 
the algorithm cannot expected to perform better than 
the usual algorithm for calculating $ \dim D(\lambda) $, via the $ \dim S(\lambda)  \times 
\dim S(\lambda)  $ matrix associated with the bilinear form on $ S(\lambda)$ -- since, after all, $ n! $ in general is 
much bigger than
$ \dim S(\lambda) $.

\medskip
\noindent
{\bf Remark}. As was pointed out to us by A. Mathas, a generator for $ D(\lambda) $ is 
given on page 41 in [J]. In our terminology it is 
$x_{\lambda} s_{\lambda} y_{\lambda^{\prime}}  s_{\lambda}^{-1}   x_{\lambda} $ and hence, by the 
arguments of the Theorem, it coincides with our generator. Our final expression of it is 
somewhat shorter, but still does not permit calculations much beyond the ones already indicated.

\medskip
\noindent
{\bf Remark}. It is known from [Mu92] that $ \coeff_1(E_{\lambda}) = \frac{1}{h_{\lambda}} $ 
where $ h_{\lambda} $ is the hook-product as above. In fact since $ R S_n $ is a symmetric algebra 
we get from Proposition 9.17 of [CR] 
that this fact also holds for $ E_t $ when $ t \in \Std(\lambda) $. 
Based on GAP calculations we conjecture that 
the coefficient of any $ w \in S_n $ in $ E_{\lambda} $ is either zero or 
on the form $ \frac{1}{\,k_w h_{\lambda}\,}$ for some nonzero integer $ k_w$. 
According to our GAP-calculations, a
similar 
statement does not hold for the general $ E_t$.


\begin{thebibliography}{99}
\bibitem[A1]. H. H. Andersen, The strong linkage principle, 
J. Reine Ang. Math. {\bf 315} (1980), 53--59.
\bibitem[A2]. H. H. Andersen, 
The Frobenius morphism on the cohomology of homogeneous vector bundles on $G/B$, 
Annals Math. {\bf 112} (1980), 113--121.
\bibitem[AMR].
S. Ariki, A. Mathas, H. Rui, Cyclotomic Nazarov–Wenzl algebras, Nagoya J. Math., {\bf 182} (2006),
47--134. (Special issue in honour of George Lusztig).
\bibitem[BK].
J. Brundan, A. Kleshchev, 
Blocks of cyclotomic Hecke algebras and Khovanov-Lauda algebras.  
Invent. Math. {\bf 178}  (2009),  no. 3.
\bibitem[BKW].
J. Brundan, A. Kleshchev,  W. Wang, 
Graded Specht modules, J. Reine Angew. Math. {\bf 655} (2011), 61--87.
\bibitem[CR]. C. W. Curtis, I. Reiner, Methods of representation theory, 
vol. II, Pure and Applied Mathematics, John Wiley, New York, 1987.
\bibitem[DJ]. R. Dipper, G. James, Representations of Hecke algebras of general linear groups. 
Proc. London Math. Soc. (3) {\bf 52} (1986), no. 1, 20--52. 
\bibitem[GL]. J. Graham, G.I. Lehrer, Cellular algebras, Inventiones Math. {\bf 123} (1996), 1--34.
\bibitem[HuMa].  J. Hu, A. Mathas,
Graded cellular bases for the cyclotomic Khovanov-Lauda-Rouquier algebras of type $A$,
Adv. Math., {\bf 225} (2010), 598--642. 
\bibitem[HuMa1]. 
J. Hu, A. Mathas, Graded induction for Specht modules, Int. Math. Res. Notices {\bf 2012}  (2012)
doi: 10.1093/imrn/rnr058.
\bibitem[J]{james} G. D. James, The representation theory of the 
symmetric groups {\bf 682}, \,Lecture notes in mathematics, Springer Verlag
(1978).
\bibitem[Jan]. J. C. Jantzen, 
Representations of Algebraic Groups, Academic Press, 2003
ISBN 13: 9780821835272.
\bibitem[Ju1]. A. A. Jucys, On the Young operators of symmetric groups, Liroosk. Fiz Sb. {\bf 6} (1966)
163--180.
\bibitem[Ju2]. 
A. A. Jucys, Factorisation of Young's projection operators of symmetric groups, Litousk.
Fiz. Sb. {\bf 11} (1971), l--10.
\bibitem[Ju3].
A. A. Jucys, Symmetric polynomials and the centre of the symmetric group ring,
Rep. Mat. Phm. { \bf 5} (1974), 107--112.
\bibitem[LLT]{lascoux leclerc thibon} A. Lascoux, B. Leclerc, J.-Y. Thibon,
Hecke algebras at roots of unity and crystal bases of quantum affine 
algebras, Commun. Math. Phys. {\bf 181} (1996), 205--263.
\bibitem[Ma]. A. Mathas, 
Iwahori-Hecke Algebras and Schur Algebras of the Symmetric Group. Univ. Lecture
Series 15, Amer. Math. Soc., 1999.
\bibitem[MaSo]. A. Mathas, Seminormal forms and Gram determinants for cellular algebras, J.
Reine Angew. Math., {\bf 619} (2008), 141-173.  With an appendix by M. Soriano.
\bibitem[Mu81].G. E. Murphy, A New Construction of Young's Seminormal 
Representation of the Symmetric Groups, Journal of Algebra {\bf 69} (1981), 
287-297.
\bibitem[Mu83].G. E. Murphy, The Idempotents of the Symmetric Groups and
Nakayama's Conjecture, Journal of Algebra {\bf 81} (1983), 258--265.
\bibitem[Mu92].G. E. Murphy, On the Representation Theory of the 
Symmetric
Groups and Associated Hecke Algebras, Journal of Algebra {\bf 152} (1992),
492-513.
\bibitem[Mu95].G. E. Murphy, The Representations of Hecke Algebras of type
$ A_n $, Journal of Algebra {\bf 173} (1995), 97--121. 
\bibitem[OV].A. Y. Okounkov, A. M. Vershik, 
A new approach to representation theory of symmetric
groups, Selecta Math., New Series {\bf 2} (1996), no. 4, 581--605.
\bibitem[RH1]. S. Ryom-Hansen, 
Grading the translation functors in type $A$, J. Algebra {\bf 274} (2004), no. 1, 138--163.
\bibitem[RH2]. S. Ryom-Hansen, 
On the denominators of Young's seminormal
basis, arXiv:0904.4243.	
\bibitem[RH3]. S. Ryom-Hansen, A $q$-analogue of Kempf's venishing Theorem, 
Mosc. Math. J., {\bf 3} (2003), no. 1, 173--187. 
\bibitem[S97]. M. Sch\"onert et~al.
  \newblock { {GAP} --
            {Groups}, {Algorithms}, and {Programming} --
            version 3 release 4 patchlevel 4"}.
  \newblock Lehrstuhl D f\"ur Mathematik,
            Rheinisch Westf\"alische
            Technische Hochschule, Aachen, Germany, 1997.
\bibitem[W]. M. Wildon, Notes on Murphy operators and Nakayamas's conjecture, available on the site
http://www.ma.rhul.ac.uk/$\sim$uvah099/Maths/Murphy.pdf



\end{thebibliography}
\end{document}